\documentclass[10pt]{article}
\usepackage{upgreek}
\usepackage{latexsym,epsfig,bm,amssymb}
\usepackage{color}
\usepackage{amsthm,mathrsfs}
\usepackage{times}

\def\version{November 10, 2007}

\DeclareSymbolFont{AMSb}{U}{msb}{m}{n}
\DeclareSymbolFontAlphabet{\mathbb}{AMSb}

\newcommand{\beqn}{\begin{eqnarray}}
\newcommand{\eeqn}{\end{eqnarray}}
\newcommand{\be}{\begin{equation}}
\newcommand{\ee}{\end{equation}}
\newcommand{\ba}{\begin{array}}
\newcommand{\ea}{\end{array}}

\newcommand{\E}{\mathscr{E}}

\newcommand{\bS}{{\bf S}}

\newcommand{\bd}{\begin{definition}}
 \newcommand{\ed}{\end{definition}}
\newcommand{\bt}{\begin{theorem}}
\newcommand{\et}{\end{theorem}}
\newcommand{\bqt}{\begin{qtheorem}}
\newcommand{\eqt}{\end{qtheorem}}

\newcommand{\bp}{\begin{pro}}
\newcommand{\ep}{\end{pro}}

\newcommand{\bl}{\begin{lemma}}
\newcommand{\el}{\end{lemma}}
\newcommand{\bc}{\begin{cor}}
\newcommand{\ec}{\end{cor}}

\newcommand{\bex}{\begin{example}}
\newcommand{\eex}{\end{example}}
\newcommand{\bexs}{\begin{examples}}
\newcommand{\eexs}{\end{examples}}

\newcommand{\bexe}{\begin{exercice}}
\newcommand{\eexe}{\end{exercice}}

\newcommand{\br}{\begin{remark} }
\newcommand{\er}{\end{remark}}

\def\xx{{m}}
\def\yy{{n}}

\newcommand\dist{\mathop{\rm dist}}
\newcommand\supp{\mathop{\rm supp}}
\newcommand\p{\partial}
\newcommand{\at}[1]{\vert\sb{\sb{#1}}}

\def\Re{{\rm Re\, }}
\def\Im{{\rm Im\,}}
\providecommand\C{{\mathbb C}}
\renewcommand\C{{\mathbb C}}
\newcommand{\R}{{\mathbb R}}
\newcommand{\N}{{\mathbb N}}

\newcommand{\Abs}[1]{\left\vert#1\right\vert}
\newcommand{\abs}[1]{\vert #1\vert}
\newcommand{\Norm}[1]{\left\Vert #1\right\Vert}
\newcommand{\norm}[1]{\Vert #1\Vert}
\newcommand{\const}{{\rm const}}

\newcommand\sothat{{\rm :}\ }
\newcommand\sgn{\mathop{\rm sgn}}

\newcommand\mod{\mathop{\rm mod}}

\providecommand{\ltor}[1]{
\ifnum #1=1{\it i}\else\ifnum #1=2{\it ii}\else\ifnum #1=3{\it iii}
\else\ifnum #1=4 {\it iv}\fi\fi\fi\fi
}

\DeclareMathSymbol{\varDelta}{\mathord}{letters}{"01}
\DeclareMathSymbol{\varTheta}{\mathord}{letters}{"02}
\DeclareMathSymbol{\varLambda}{\mathord}{letters}{"03}
\DeclareMathSymbol{\varPhi}{\mathord}{letters}{"08}
\DeclareMathSymbol{\varPsi}{\mathord}{letters}{"09}
\DeclareMathSymbol{\varOmega}{\mathord}{letters}{"0A}

\theoremstyle{plain}
\newtheorem{theorem}{Theorem}[section]

\newtheorem{lemma}[theorem]{Lemma}
\newtheorem{corollary}[theorem]{Corollary}
\newtheorem{proposition}[theorem]{Proposition}

\theoremstyle{definition}
\newtheorem{definition}[theorem]{Definition}
\newtheorem{assumption}{Assumption}

\theoremstyle{remark}
\newtheorem{remark}[theorem]{Remark}
\newtheorem{example}[theorem]{Example}

\makeatletter\@addtoreset{equation}{section}
\makeatletter\@addtoreset{theorem}{section}
\makeatother

\begin{document}


\title{
On global attraction to solitary waves for the
Klein-Gordon field coupled to several nonlinear oscillators
}


\author{
{\sc Alexander Komech}
\footnote{On leave from 
Institute of the Information Transmission Problems, RAS.
Supported in part by Alexander von Humboldt Research Award,
by Max-Planck Institute for Mathematics in the Sciences (Leipzig),
and by grants FWF P19138-N13,
DFG 436 RUS 113/929/0-1, and RFBR 07-01-00018a.}
\\
{\it\small Faculty of Mathematics, University of Vienna, Wien A-1090, Austria}
\\
{\it\small
Institute for Information Transmission Problems,
Moscow 101447, Russia}
\\ \\
{\sc Andrew Komech}
\footnote{Supported in part
by Max-Planck Institute for Mathematics in the Sciences (Leipzig) and
by the National Science Foundation
under Grant DMS-0600863.
}
\\
{\it\small
Mathematics Department, Texas A\&M University,
College Station, TX 77843, USA}
\\
{\it\small
Institute for Information Transmission Problems,
Moscow 101447, Russia}
}
\date{\version}

\maketitle
\begin{abstract}
The global attraction is established for all finite energy solutions
to a model $\mathbf{U}(1)$-invariant nonlinear Klein-Gordon equation
in one dimension coupled to a finite number of nonlinear oscillators:
We prove that {\it each finite energy solution}
converges as $t\to\pm\infty$ to the set of all ``nonlinear eigenfunctions''
of the form $\phi(x)e\sp{-i\omega t}$
if  all oscillators are strictly nonlinear, and
the distances between neighboring oscillators are sufficiently small.


Our approach is based on the analysis of {\it omega-limit trajectories}
which form the global attractor.
We show that their time spectrum is a priori compact.
Then
the nonlinear spectral analysis based on the Titchmarsh convolution theorem
allows to reduce the time-spectrum to one point.
This implies that each omega-limit trajectory is a solitary wave.
Physically, the {\it global attraction} is caused by the nonlinear
energy transfer from lower harmonics to the continuous spectrum
and subsequent dispersive radiation.
The Titchmarsh theorem allows to prove that this energy transfer
and radiation are absent only for the solitary waves.

To check the sharpness of our conditions,
we construct counterexamples showing
the global attractor can contain
``multifrequency solitary waves''
if the distance between oscillators is large 
or if some of them are linear.
\end{abstract}


\section{Introduction}
The long time asymptotics for
nonlinear
wave equations
have been the subject of intensive research,
starting with the pioneering papers by
Segal \cite{MR0153967,MR0152908},
Strauss \cite{MR0233062},
and Morawetz and Strauss \cite{MR0303097},
where the
nonlinear scattering and the local attraction to zero solution
were proved.
Local attraction to solitary waves,
or \emph{asymptotic stability},
in
$\mathbf{U}(1)$-invariant dispersive systems
was addressed in
\cite{MR1071238,MR1199635e,MR1170476,MR1334139}
and then developed in
\cite{MR1488355,MR1681113,MR1893394,MR1835384,MR1972870,MR2027616}.
Global attraction
to \emph{static},
stationary solutions in the dispersive systems
\emph{without $\mathbf{U}(1)$ symmetry}
was established in
\cite{MR1203302e,MR1359949,MR1412428,MR1434147,MR1726676,MR1748357}.

We would like to have
the dynamical description
of the Bohr transitions to quantum stationary states
in coupled nonlinear systems
of Quantum Physics.
This suggests investigation of the global attractors in
nonlinear Hamiltonian hyperbolic equations
with $\mathbf{U}(1)$-symmetry
(see \cite{ubk-arma} for the discussion).
The first result about
the global attraction to solitary waves in
a model with these properties
was obtained in \cite{ubk-cr,ubk-arma},
where we considered the Klein-Gordon equation coupled to
one nonlinear oscillator.

We are aware of only one other recent advance \cite{MR2304091}
in the field
of nonzero global attractors for Hamiltonian PDEs.
In that paper, the global attraction for the nonlinear Schr\"odinger equation
in dimensions $n\ge 5$ was considered.
The dispersive wave was explicitly specified
using the rapid decay of local energy in higher dimensions.
The global attractor was proved to be compact, but it was
neither identified with the set of solitary waves nor was proved
to be of finite dimension \cite[Remark 1.18]{MR2304091}.

In the present paper, we prove the attraction
to the set of solitary waves
for all finite energy solutions to the Klein-Gordon equation
coupled to any finite number of nonlinear oscillators.
For the proof, we develop an approach of the spectral inflation
\cite{ubk-arma}
justified by the Titchmarsh Convolution Theorem.
This justification requires new arguments and appropriate conditions.
We demonstrate the sharpness of these conditions constructing
counterexamples.

Our model is based on the complex Klein-Gordon field $\psi(x,t)$, 
interacting with $N$ nonlinear oscillators
located at the points
$X\sb{1}<X\sb{2}<\dots<X\sb{N}$:
\begin{equation}\label{kg-no-0}
\ddot\psi
=\psi''-m^2\psi+\sum\sb{J}\delta(x-X\sb{J})F\sb{J}(\psi(X\sb{J},t)),
\qquad x\in\R,
\end{equation}
where
$m>0$
and
$F\sb{J}$ are nonlinear functions
describing nonlinear oscillators at the points $X\sb{J}$.
The dots stand for the derivatives in $t$,
and the primes for the derivatives in $x$.
All derivatives and the equation are understood in
the sense of distributions.
We assume that equation (\ref{kg-no-0}) is $\mathbf{U}(1)$-invariant;
that is,
\begin{equation}\label{inv-f}
F\sb{J}(e\sp{i\theta}\psi)=e\sp{i\theta} F\sb{J}(\psi),
\qquad\theta\in\R,\quad\psi\in\C,\quad 1\le J\le N.
\end{equation}
This symmetry leads to the charge conservation
and to the existence of the solitary wave solutions,
which are finite energy solutions of the following form:
\begin{equation}\label{soliton}
\psi\sb\omega(x,t)=\phi\sb\omega(x)e\sp{-i\omega t},\quad\omega\in\R,
\quad\phi\sb\omega\in H\sp{1}(\R).
\end{equation}
Above, $H\sp{1}(\R)$ denotes the Sobolev space.

\begin{definition}\label{dSS}
$S$
is the set of all functions
$\phi\sb\omega(x)\in H\sp{1}(\R)$ with $\omega\in\R$,
so that $\phi\sb\omega(x)e\sp{-i\omega t}$
is a solution to (\ref{kg-no-0}).
\end{definition}
Note that $S$ also contains the zero solution.

Generically,
the factor-space $S/\mathbf{U}(1)$
is isomorphic to a finite union of one-dimensional
intervals.
The set of all solitary waves for equation (\ref{kg-no-0})
is described in Proposition~\ref{prop-solitons}.
Typically, such solutions exist for
$\omega$ from an interval or a collection of intervals
of the real line.

\bigskip

Our main result is the following long-time asymptotics:
In the case when
all oscillators are polynomial and strictly nonlinear
(see Assumptions~\ref{ass-a} and \ref{ass-nonl} below)
and all distances $\abs{X\sb{J+1}-X\sb{J}}$ are sufficiently small,
we prove
that any finite energy solution converges
to the set $S$ of all solitary waves:
\begin{equation}\label{attraction}
\psi(\cdot,t)\longrightarrow S,
\qquad t\to\pm\infty,
\end{equation}
where the convergence holds in local energy seminorms.

Let us give a brief sketch of our approach.
We introduce a concept of the omega-limit trajectories
$\beta(x,t)$ which 
play a crucial role in the proof.
We define omega-limit trajectories as the limits
\[
\psi(x,t+s\sb j)\rightarrow\beta(x,t),
\qquad
(x,t)\in\R^2,
\]
for some sequence of times
$s\sb j\to+\infty$.
We will prove that all omega-limit trajectories
are solitary waves, thus finishing the proof.
To complete this program,
we study the time spectrum of solutions,
that is, their Fourier-Laplace transform in time.
We need to prove that
$\beta(x,t)=\phi\sb\omega(x)e^{-i\omega t}$, that is,
that the time spectrum of $\beta$ consists of at most one frequency.
First, we show that the spectrum of the solution
at $x=X\sb{1}$ and $x=X\sb{N}$
is absolutely continuous for $\abs{\omega}>m$.
At the points $x\in(X\sb 1,X\sb N)$,
the nonlinearity may extend the singular part of the spectrum
to be at most $[-\Lambda,\Lambda]$,
for some bounded $\Lambda$.
Outside of this interval, the spectrum is absolutely continuous.
This allows to prove that the spectrum
of any omega-limit trajectory at $x=X\sb{1}$ and $x=X\sb{N}$
is contained in $[-m,m]$,
while at the points $x\in(X\sb 1,X\sb N)$
the spectrum is contained in $[-\Lambda,\Lambda]$.
The next important observation is
that each omega-limit trajectory
is also a solution to the original nonlinear Klein-Gordon equation.
This allows to apply the Titchmarsh theorem and prove
that the spectrum of any omega-limit trajectory
at all points $x\in\R$ consists of at most one frequency.
At this last step,
one needs
the assumptions that the oscillators are strictly nonlinear
and located sufficiently close to one another.

The requirement that the nonlinearities $F\sb{J}$ are polynomial
allows us to apply the Titchmarsh theorem
which is vital in the proof.
We construct counterexamples showing the sharpness of our assumptions
for the global attraction to the solitary waves. 
Namely, for $N=2$,
we construct multifrequency solitary waves
in the case when the distance $\abs{X\sb 2-X\sb 1}$ is sufficiently large
or one of the oscillators is linear.

Let us mention that
in the case of $N$ oscillators,
considered in this paper,
the general plan of the proof is similar
to the case of one oscillator (see \cite{ubk-cr,ubk-arma}).
However, the justifications of all steps 
are based on new arguments.
In particular, the application of the Titchmarsh theorem
required a new construction.

\bigskip

Our paper is organized as follows.
In Section~\ref{sect-results}, we formulate our main results.
In Section~\ref{sect-splitting}, we separate the first dispersive component. 
In Sections~\ref{sect-spec} and ~\ref{sect-ac}, 
we construct spectral representation 
for the remaining component, 
and prove absolute continuity of its spectrum for high frequencies.
In Sections~\ref{sect-bound}, we  
separate the second dispersive component corresponding to the 
high frequencies
and establish compactness for the remaining \emph{bound} component
with the bounded spectrum.
In Section~\ref{sect-nonlinear-analysis},
we study omega-limit trajectories of the solution.
In Section~\ref{sect-counterexamples} we collect counterexamples,
and in Appendix~\ref{sect-existence}
we establish global well-posedness.

\section{Main results}
\label{sect-results}

\subsection*{Model}
We consider the Cauchy problem for the Klein-Gordon equation
with the nonlinearity concentrated at the points
$X\sb{1}<X\sb{2}<\dots<X\sb{N}$:
\begin{equation}\label{kg-no}
\left\{
\begin{array}{l}
\ddot\psi(x,t)
=\psi''(x,t)-m^2\psi(x,t)
+\sum\sb{J}\delta(x-X\sb{J})F\sb{J}(\psi(X\sb{J},t)),
\qquad
x\in\R,
\quad t\in\R,
\\
\psi\at{t=0}=\psi\sb{0}(x),
\qquad
\dot\psi\at{t=0}=\pi\sb{0}(x).
\end{array}\right.
\end{equation}
If we identify a complex number $\psi=u+i v\in\C$
with the two-dimensional vector
$(u,v)\in\R\sp 2$,
then, physically, equation (\ref{kg-no}) describes small crosswise
oscillations of the infinite
string in three-dimensional space
$(x,u,v)$
stretched along the $x$-axis.
The string is subject to
the action of
an ``elastic force'' $-m^2\psi(x,t)$ and
coupled to nonlinear oscillators
of forces $F\sb{J}(\psi)$
attached at the points $X\sb{J}$.
We denote by $\mathcal{X}$
the set of all the locations of oscillators:
\begin{equation}\label{def-x}
\mathcal{X}=\{X\sb 1,\,X\sb 2,\,\dots,\,X\sb{N}\}.
\end{equation}
We will assume that
the oscillator forces $F\sb{J}$ admit real-valued potentials:
\begin{equation}\label{P}
F\sb{J}(\psi)=-\nabla U\sb{J}(\psi),\quad\psi\in\C,
\qquad
U\sb{J}\in C\sp 2(\C),
\end{equation}
where the gradient is taken with respect to $\Re\psi$ and $\Im\psi$.
We define
$\Psi(t)=
\left[\!\scriptsize{\begin{array}{c}
\psi(x,t)\\\pi(x,t)\end{array}}\!\right]$
and write the Cauchy problem
(\ref{kg-no})
in the vector form:
\begin{equation}\label{kg-no-cp}
\dot\Psi(t)
=
\left[\begin{array}{cc}0&1\\\p\sb x\sp 2-m^2&0\end{array}\right]
\Psi(t)
+
\sum\sb{J}
\delta(x-X\sb{J})\left[\begin{array}{c}0\\F\sb{J}(\psi)\end{array}\right],
\qquad
\Psi\at{t=0}
=\Psi\sb{0}
\equiv\left[\begin{array}{c}\psi\sb{0}\\\pi\sb{0}\end{array}\right].
\end{equation}
Equation (\ref{kg-no-cp})
formally can be written as a Hamiltonian system,
\begin{equation}\label{kg-no-h}
\dot\Psi(t)=\mathcal{J}\,D\mathcal{H}(\Psi),
\qquad
\mathcal{J}=\left[\begin{array}{cc}0&1\\-1&0\end{array}\right],
\end{equation}
where $D\mathcal{H}$ is the variational
derivative of the Hamilton functional
\begin{equation}\label{hamiltonian}
\mathcal{H}(\Psi)
=\frac 1 2
\int\limits\sb{\R}
\left(
\abs{\pi}\sp 2+\abs{\psi'}\sp 2+m^2\abs{\psi}\sp 2
\right)
dx
+\sum\sb{J} U\sb{J}(\psi(X\sb{J})),
\qquad
\Psi=\left[\begin{array}{c}\psi(x)\\\pi(x)\end{array}\right].
\end{equation}
We
assume that the potentials $U\sb{J}(\psi)$ are $\mathbf{U}(1)$-invariant,
where $\mathbf{U}(1)$ stands for the unitary group
$e\sp{i\theta}$, $\theta\in\R\mod 2\pi$.
Namely, we assume that
there exist $u\sb{J}\in C\sp 2(\R)$ such that
\begin{equation}\label{inv-u}
U\sb{J}(\psi)=u\sb{J}(\abs{\psi}\sp 2),
\qquad\psi\in\C,
\quad 1\le J\le N.
\end{equation}

\begin{remark}
In the context of the model of the infinite string in $\R\sp{3}$
that we described after (\ref{kg-no}),
the assumption (\ref{inv-u})
means that the potentials $U\sb{J}(\psi)$
are rotation-invariant with respect to the $x$-axis.
\end{remark}

Conditions (\ref{P}) and (\ref{inv-u})
imply that
\begin{equation}\label{def-a}
F\sb{J}(\psi)=
\alpha\sb{J}(\abs{\psi}^2)\psi,
\qquad\psi\in\C,
\end{equation}
where
$\alpha\sb{J}(\cdot)
=-2 u\sb{J}'(\cdot)\in C\sp 1(\R)$
are real-valued.
Therefore, (\ref{inv-f}) holds.
Since (\ref{kg-no-cp}) is $\mathbf{U}(1)$-invariant,
the N\"other theorem formally implies that the {\it charge functional}
\begin{equation}\label{cal-Q}
\mathcal{Q}(\Psi)
=\frac{i}{2}\int\sb{\R}
\left(\overline\psi\pi-\overline{\pi}\psi\right)\,dx,
\qquad
\Psi=\left[\begin{array}{c}\psi(x)\\\pi(x)\end{array}\right],
\end{equation}
is conserved for solutions $\Psi(t)$ to (\ref{kg-no-cp}).

Let us introduce
the phase space ${\E}$
of finite energy states for equation (\ref{kg-no}).
Denote by $L\sp 2$ the complex Hilbert space $L\sp 2(\R)$
with the norm $\norm{\cdot}\sb{L\sp 2}$,
and denote by $\norm{\cdot}\sb{L\sp 2\sb R}$ the norm in $L\sp 2(-R,R)$
for $R>0$.

\begin{definition}
\begin{enumerate}
\item
${\E}$ is the Hilbert space of the states
$\Psi=(\psi,\pi)$,
with the norm
\begin{equation}\label{def-e}
\norm{\Psi}\sb{\E}^2
:=
\norm{\pi}\sb{L\sp 2}^2
+\norm{\psi'}\sb{L\sp 2}^2+m^2\norm{\psi}\sb{L\sp 2}^2.
\end{equation}
\item
${\E}\sb{F}$ is the space $\E$
endowed with the Fr\'echet topology
defined by local energy seminorms
\begin{equation}\label{def-e-r}
\norm{\Psi}\sb{\E,R}^2
:=
\norm{\pi}\sb{L\sp 2(-R,R)}^2
+\norm{\psi'}\sb{L\sp 2(-R,R)}^2+m^2\norm{\psi}\sb{L\sp 2(-R,R)}^2,
\qquad
R>0.
\end{equation}
\end{enumerate}
\end{definition}

\begin{remark}
The space $\E\sb{F}$ is metrizable.
The metric could be introduced by
\begin{equation}\label{def-e-metric}
\dist(\Psi,\Phi)
=\sum\sb{R=1}\sp{\infty} 2^{-R}
\norm{\Psi-\Phi}\sb{\E,R}.
\end{equation}
\end{remark}

Equation (\ref{kg-no-cp})
is formally a Hamiltonian system with
the phase space ${\E}$
and the Hamilton functional $\mathcal{H}$.
Both
$\mathcal{H}$
and $\mathcal{Q}$ are continuous functionals on ${\E}$.
Let us note that
${\E}={H\sp{1}}\oplus L\sp 2$,
where $H\sp{1}$ denotes the Sobolev space
\[
H\sp{1}=H\sp{1}(\R)
=\{\psi(x)\in L\sp 2(\R):\;\psi'(x)\in L\sp 2(\R)\}.
\]
We introduced into (\ref{def-e})
the factor $m^2>0$,
to have a convenient relation
$\mathcal{H}(\psi,\dot\psi)
=\frac 1 2\norm{(\psi,\dot\psi)}\sb{\E}^2+\sum\sb{J} U\sb{J}(\psi(X\sb{J}))$.

\subsection*{Global well-posedness}

To have a priori estimates available for the proof of the global
well-posedness, we assume that
\begin{equation}\label{bound-below}
U\sb{J}(\psi)\ge {A}\sb{J}-{B}\sb{J}\abs{\psi}^2
\quad{\rm for}\ \psi\in\C,\quad
{\rm where}\quad {A}\sb{J}\in\R,\quad {B}\sb{J}\ge 0,
\quad
1\le J\le N;
\quad
\sum\sb{J} {B}\sb{J}<m.
\end{equation}

\begin{theorem}\label{theorem-well-posedness}
Let $F\sb{J}(\psi)$ satisfy conditions (\ref{P})
and (\ref{inv-u}):
\[
F\sb{J}(\psi)=-\nabla U\sb{J}(\psi),\qquad
U\sb{J}(\psi)=u\sb{J}(\abs{\psi}^2),
\qquad u\sb{J}(\cdot)\in C\sp 2(\R).
\]
Additionally,
assume that (\ref{bound-below}) holds.
Then:
\begin{enumerate}
\item
For every $\Psi\sb{0}\in {\E}$ the Cauchy problem
(\ref{kg-no-cp}) has a unique solution
$\Psi(t)$ such that $\Psi\in C(\R,{\E})$.
\item
The map
$W(t):\;\Psi\sb{0}\mapsto\Psi(t)$
is continuous in ${\E}$
for each $t\in\R$.
\item
The energy and charge are conserved:
$\mathcal{H}(\Psi(t))=\const$,
$\mathcal{Q}(\Psi(t))=\const$,
$t\in\R$.
\item
The following \emph{a priori} bound holds:
$
\norm{\Psi(t)}\sb{\E}
\le C(\Psi\sb{0})$,
$t\in\R$.
\end{enumerate}
\end{theorem}

We prove this Theorem in Appendix~\ref{sect-existence}.

\subsection*{Solitary waves and the main theorem}

\begin{definition}\label{def-solitary-waves}
\begin{enumerate}
\item
The solitary waves of equation (\ref{kg-no})
are solutions of the form
\begin{equation}\label{solitary-waves}
\psi(x,t)=\phi\sb\omega(x)e\sp{-i\omega t},
\qquad
{\rm where}
\quad
\omega\in\R,
\quad
\phi\sb\omega\in H\sp{1}(\R).
\end{equation}
\item
The solitary manifold
is the set
$
\bS
=
\left\{
(\phi\sb\omega,-i\omega\phi\sb\omega)
\sothat\omega\in\R,\ \phi\sb\omega\in H\sp 1(\R)
\right\}
\subset\E.
$
\end{enumerate}
\end{definition}

\begin{remark}
\begin{enumerate}
\item
Identity (\ref{inv-f}) implies that the set
$\bS$
is invariant under multiplication by $e\sp{i\theta}$,
$\theta\in\R$.
\item
Let us note that for any $\omega\in\R$
there is a zero solitary wave with
$\phi\sb\omega(x)\equiv 0$
since $F\sb{J}(0)=0$ by (\ref{def-a}).
\item
According to (\ref{def-a}),
$\alpha\sb{J}(\abs{C}^2)=F\sb{J}(C)/C\in\R$
for any $C\in\C\backslash 0$.
\end{enumerate}
\end{remark}

\begin{definition}\label{def-sn}
The function $F\sb{J}(\psi)$ is \emph{strictly nonlinear}
if the equation
$\alpha\sb{J}(C^2)=a$
has a discrete
(or empty)
set of positive roots
$C$ for each particular $a\in\R$.
\end{definition}

The following proposition provides a concise description
of all solitary waves.
Formally this proposition 
is not necessary for our exposition.
\begin{proposition}
\label{prop-solitons}
Assume that $F\sb{J}(\psi)$ satisfy (\ref{inv-f}) and that
$F\sb{J}(\psi)$, $1\le J\le N$, are strictly nonlinear
in the sense of Definition~\ref{def-sn}.
Then all solitary wave solutions
to {(\ref{kg-no})}
are given by {(\ref{solitary-waves})} with
\begin{equation}\label{solitary-wave-profile}
\phi\sb\omega(x)
=\sum\sb{J}
C\sb{J} e^{-\kappa(\omega)\abs{x-X\sb{J}}},
\qquad
\kappa(\omega)=\sqrt{m^2-\omega^2},
\end{equation}
where
$\omega\in[-m,m]$
and $C\sb{J}\in\C$,
$1\le J\le N$,
satisfy the following relations:
\begin{equation}\label{kaka}
2\kappa(\omega) C\sb{J}
=F\sb{J}\Big(\sum\sb{K}C\sb K e^{-\kappa(\omega)\abs{X\sb{J}-X\sb K}}\Big).
\end{equation}
\end{proposition}

\begin{remark}\label{remark-zero}
By (\ref{solitary-wave-profile}),
$\omega=\pm m$ can only correspond to zero solution.
\end{remark}

The proof of this Proposition repeats the proof
of a similar result for the case $N=1$
in \cite{ubk-arma}.

As we mentioned before,
we need to assume that the nonlinearities are
nonlinear polynomials.
This condition
is crucial in our argument:
It will allow to apply the Titchmarsh convolution theorem.

\bigskip

Let us formulate all the assumptions which we need
to formulate the main result.

\begin{assumption}\label{ass-a}
For all $1\le J\le N$,
\begin{equation}\label{f-is-such}
F\sb{J}(\psi)=-\nabla U\sb{J}(\psi),
\quad{\rm where}\quad
U\sb{J}(\psi)=\sum\limits\sb{n=0}\sp{p\sb{J}}u\sb{J,n}
\abs{\psi}\sp{2n}~,
\qquad
\ \ u\sb{J,n}\in\R.
\end{equation}
\end{assumption}

\begin{assumption}\label{ass-nonl}
For all $1\le J\le N$, we have
\begin{equation}\label{nonli}
u\sb{J,p\sb{J}}>0
\quad
{\rm and}
\quad
p\sb{J}\ge 2.
\end{equation}
\end{assumption}

Assumptions~\ref{ass-a} and \ref{ass-nonl} guarantee that
all nonlinearities $F\sb{J}$ are strictly nonlinear
and satisfy (\ref{P}), (\ref{inv-u}),
and also that the bound (\ref{bound-below}) takes place.

We introduce the following quantities:
\begin{equation}\label{def-mu}
\mu\sb{1}=m,
\quad
\mu\sb{J+1}=(2p\sb{J}-1)\mu\sb{J};
\qquad
\mu'\sb{N}=m,
\quad
\mu'\sb{J}=(2p\sb{J+1}-1)\mu'\sb{J+1},
\qquad
1\le J\le N-1,
\end{equation}
where $p\sb J$ are exponentials from (\ref{f-is-such}).
We also denote
\begin{equation}\label{def-m-lambda}
\Lambda=\max\limits\sb{1\le J\le N}(2p\sb{J}-1){M}\sb{J},
\qquad
{\rm where}
\quad
{M}\sb{J}=\min(\mu\sb{J},\mu'\sb{J}).
\end{equation}


\begin{assumption}\label{ass-delta-small}
The intervals
$[X\sb{J},X\sb{J+1}]$, $\ 1\le J\le N-1$,
are small enough so that
\begin{equation}\label{delta-small}
\Lambda<
\sqrt{\frac{\pi^2}{\abs{X\sb{J+1}-X\sb{J}}^2}+m^2},
\qquad
1\le J\le N-1.
\end{equation}
\end{assumption}

Our main result is the following theorem.

\begin{theorem}[Main Theorem]
\label{main-theorem}
Let Assumptions~\ref{ass-a}, ~\ref{ass-nonl}, and
\ref{ass-delta-small} hold.
Then for any $\Psi\sb{0}\in\E$
the solution $\Psi(t)\in C(\R,\E)$
to the Cauchy problem {(\ref{kg-no-cp})}
converges to
$\bS$:
\begin{equation}\label{cal-A}
\lim\sb{t\to\pm\infty}\dist(\Psi(t),\bS)=0,
\end{equation}
where
$\dist(\Psi,\bS):=\inf\limits\sb{\bm s\in\bS}\dist(\Psi,\bm s)$,
and $\dist$ is introduced in (\ref{def-e-metric}).
\end{theorem}

\begin{remark}
\begin{enumerate}
\item
The solution $\Psi(t)$ exists by Theorem~\ref{theorem-well-posedness}
since Assumptions~\ref{ass-a} and \ref{ass-nonl}
guarantee that conditions 
(\ref{P}), (\ref{inv-u}), and (\ref{bound-below})
hold.
\item
It suffices to prove Theorem~\ref{main-theorem}
for $t\to +\infty$.
\item
In Sections~\ref{sect-example-w} and \ref{sect-example-li},
we construct counterexamples to the
convergence (\ref{cal-A})
in the case when 
Assumption~\ref{ass-nonl}
or Assumption~\ref{ass-delta-small}
are not satisfied.
\item
For the real initial data,
we obtain a real-valued solution $\psi(t)$ to (\ref{kg-no}).
Therefore, the convergence (\ref{cal-A})
of $\Psi(t)=(\psi(t),\dot\psi(t))$
to the set of pairs
$(\phi\sb\omega,-i\omega\phi\sb\omega)$
with $\omega\in\R$
implies that $\Psi(t)$
locally converges to zero:
\[
\lim\sb{t\to\infty}\dist(\Psi(t),0)=0.
\]
\end{enumerate}
\end{remark}

\section{Separation of dispersive component}
\label{sect-splitting}

Let us split the solution
$\psi(x,t)$
into two components,
$\psi(x,t)=\chi(x,t)+\varphi(x,t)$,
which are defined for all $t\in\R$
as solutions to the following Cauchy problems:
\begin{eqnarray}
&&
\ddot\chi(x,t)=\chi''(x,t)-m^2\chi(x,t),
\qquad
(\chi,\dot\chi)\at{t=0}=(\psi\sb{0}(x),\pi\sb{0}(x)),
\label{kg-no-cp-1}
\\
\nonumber
\\
&&
\ddot\varphi(x,t)=\varphi''(x,t)-m^2\varphi(x,t)
+
\sum\sb{J}
\delta(x-X\sb{J})f\sb{J}(t),
\qquad
(\varphi,\dot\varphi)\at{t=0}=(0,0),
\label{kg-no-cp-2-0}
\end{eqnarray}
where
$(\psi\sb{0}(x),\pi\sb{0}(x))$ is the initial data from (\ref{kg-no}),
and
\begin{equation}\label{def-f}
f\sb{J}(t):=F\sb{J}(\psi(X\sb{J},t)),
\qquad
t\in\R.
\end{equation}
The following lemma is proved in \cite[Lemma 3.1]{ubk-arma}.
\begin{lemma}\label{lemma-decay-psi1}
There is a local energy decay for $\chi$:
\begin{equation}\label{dp0}
\lim\sb{t\to\infty}
\Norm{(\chi(\cdot,t),\dot\chi(\cdot,t))}\sb{{\E},R}=0,
\qquad\forall R>0.
\end{equation}
\end{lemma}

Let $k(\omega)$ be the analytic function
with the domain 
$D:=\C\backslash((-\infty,-m]\cup[m,+\infty))$
such that
\begin{equation}\label{def-k}
k(\omega)=\sqrt{\omega\sp 2-m^2},
\qquad\Im k(\omega)>0,
\qquad
\omega
\in D.
\end{equation}
Let us also denote its limit values for $\omega\in\R$ by
\begin{equation}\label{def-k-plus}
k\sb{\pm}(\omega):=k(\omega\pm i0),
\qquad
\omega\in\R.
\end{equation}

\bigskip

\begin{figure}[htbp]
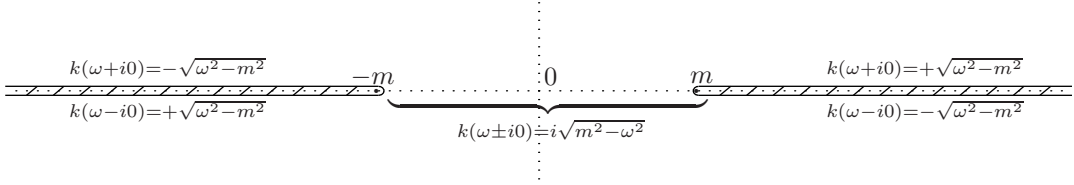

\input ukk-domain.tex
\caption{Domain $D$
and the values of $k\sb{\pm}(\omega):=k(\omega\pm i0)$, $\omega\in\R$.}
\label{fig-domain}
\end{figure}
As illustrated on Figure~\ref{fig-domain} 
(where all square roots take positive values),
we have
\begin{equation}\label{k-plus-minus}
k\sb{-}(\omega)=k\sb{+}(\omega)
\quad{\rm for}\quad
-m\le\omega\le m,
\qquad
k\sb{-}(\omega)=-k\sb{+}(\omega)
\quad{\rm for}\quad
\omega\in\R\backslash(-m,m),
\end{equation}
and also
\begin{equation}\label{omega-k}
\omega\,k\sb{+}(\omega)\ge 0
\qquad
{\rm for}
\quad
\omega\in\R\backslash(-m,m).
\end{equation}

We set $\mathcal{F}\sb{t\to\omega}[g(t)]
=\displaystyle\int\sb{\R} e^{i\omega t}g(t)\,dt$
for a function $g(t)$ from the Schwartz space $\mathscr{S}(\R)$.
Let us study the Fourier transform 
$\hat\chi(x,\omega):=\mathcal{F}\sb{t\to\omega}[\chi(x,t)]$,
which is a continuous function of $x$
valued in tempered distributions.

\begin{lemma}\label{lemma-chi-h1}
\begin{itemize}
\item
$\hat\chi(x,\omega)$
is a continuous function of 
$x\in\R$ with values in 
$L^1\sb{\rm loc}(\R)$, and 
\begin{equation}\label{nul}
\hat\chi(x,\omega)=0,
\qquad
\abs{\omega}<m.
\end{equation}
\item
The following bound holds:
\begin{equation}
\sup\limits\sb{x\in\R}~\int\sb{\abs{\omega}>m}
\abs{\hat\chi(x,\omega)}^2\omega\,k\sb{+}(\omega)\,d\omega<\infty.
\end{equation}
\end{itemize}
\end{lemma}

\begin{proof}
Set $\omega(k)=\sgn k\sqrt{m^2+k^2}$
for $k\in\R$.
Note that 
the function $k\sb{+}(\omega)$
for $\abs{\omega}>m$
is inverse to the function $\omega(k)$, $k\ne 0$.
We have:
\begin{equation}
\chi(x,t)=\frac{1}{2\pi}\int\sb{\R}
e^{-i k x}\Big[
\hat\psi\sb{0}(k)\cos(\omega(k)t)
+\hat\pi\sb{0}(k)\frac{\sin(\omega(k)t)}{\omega(k)}
\Big]\,dk.
\end{equation}
Hence, for the Fourier transform of $\chi(x,t)$,
we obtain, for any $x\in\R$:
\begin{eqnarray}
&&
\hat\chi(x,\omega)
=\int\sb{\R}
e^{-i k x}
\Big[
\hat\psi\sb{0}(k)\frac{\delta(\omega-\omega(k))+\delta(\omega+\omega(k))}{2}
+\hat\pi\sb{0}(k)\frac{\delta(\omega-\omega(k))-\delta(\omega+\omega(k))}
{2i\omega(k)}
\Big]
\,dk
\nonumber
\\
&&
\qquad
=\int\sb{\abs{\omega'}>m}
e^{-i k\sb{+}(\omega') x}
\Big[
\hat\psi\sb{0}(k\sb{+}(\omega'))\frac{\delta(\omega-\omega')+\delta(\omega+\omega')}{2}
+\hat\pi\sb{0}(k\sb{+}(\omega'))\frac{\delta(\omega-\omega')-\delta(\omega+\omega')}
{2i\omega'}
\Big]
\,\frac{\omega'\,d\omega'}{k\sb{+}(\omega')}.
\nonumber
\end{eqnarray}
The above relation is understood in the sense of distributions
of $\omega\in\R$.
We used the substitution $k=k\sb{+}(\omega')$.
Now (\ref{nul}) is obvious.
Evaluating the last integral, we get:
\[
\hat\chi(x,\omega)
=
\frac{\omega}{2k\sb{+}(\omega)}
\left\{
e^{-i k\sb{+}(\omega) x}
\hat\psi\sb{0}(k\sb{+}(\omega))
+
e^{i k\sb{+}(\omega) x}
\hat\psi\sb{0}(-k\sb{+}(\omega))
+
e^{-i k\sb{+}(\omega) x}
\frac{\hat\pi\sb{0}(k\sb{+}(\omega))}{i\omega}
-
e^{i k\sb{+}(\omega) x}
\frac{\hat\pi\sb{0}(-k\sb{+}(\omega))}{i\omega}
\right\},
\qquad\abs{\omega}>m.\nonumber
\]
We took into account that
$k\sb{+}(-\omega)=-k\sb{+}(\omega)$ for $\omega\in\R\backslash(-m,m)$
(see (\ref{k-plus-minus})).
Thus, we have:
\[
\int\limits\sb{\abs{\omega}>m}
\abs{\hat\chi(x,\omega)}^2
\omega\,k\sb{+}(\omega)\,d\omega
\le
\int\limits\sb{\abs{\omega}>m}
\Big[
\frac{\omega^2
\abs{\hat\psi\sb{0}(k\sb{+}(\omega))}^2}{k\sb{+}^2(\omega)}
+
\frac{\abs{\hat\pi\sb{0}(k\sb{+}(\omega))}^2}{k\sb{+}^2(\omega)}
\Big]
\omega\,k\sb{+}(\omega)\,d\omega
=
\int\limits\sb{\R}
\Big[
\abs{\hat\psi\sb{0}(k)}^2
+
\frac{\abs{\hat\pi\sb{0}(k)}^2}
{\omega^2(k)}
\Big]
\omega^2(k)\,dk.
\]
The finiteness of the right-hand side
follows from the finiteness of the energy of the
initial data $(\psi\sb{0},\pi\sb{0})$:
\[
\norm{(\psi\sb{0},\pi\sb{0})}\sb{\E}^2
=
\frac{1}{2\pi}
\int\sb{\R}
\Big[
\omega^2(k)
\abs{\hat\psi\sb{0}(k)}^2
+
\abs{\hat\pi\sb{0}(k)}^2
\Big]
\,dk
<\infty.
\]
\end{proof}

\section{Spectral representation}
\label{sect-spec}
The function $\varphi(x,t)=\psi(x,t)-\chi(x,t)$
satisfies the following Cauchy problem:
\begin{equation}
\ddot\varphi(x,t)=\varphi''(x,t)-m^2\varphi(x,t)
+
\sum\sb{J}
\delta(x-X\sb{J})f\sb{J}(t),
\qquad
(\varphi,\dot\varphi)\at{t=0}=(0,0),
\label{kg-no-cp-2}
\end{equation}
with $f\sb{J}(t)$ defined in (\ref{def-f}).
Note that
$\psi(X\sb{J},\cdot)\in C\sb{b}(\R)$
for $1\le J\le N$
by the Sobolev embedding,
since $(\psi(x,t),\dot\psi(x,t))\in C\sb{b}(\R,\E)$
by Theorem~\ref{theorem-well-posedness}~({\it iv}).
Hence, $f\sb J(t)\in C\sb{b}(\R)$.
On the other hand,
since $\chi(x,t)$
is a finite energy solution to the free Klein-Gordon equation,
we also have
\begin{equation}\label{psi-1-bounds}
(\chi(x,t),\dot\chi(x,t))\in C\sb{b}(\R,\E).
\end{equation}
Therefore, the function $\varphi(x,t)=\psi(x,t)-\chi(x,t)$ satisfies
\begin{equation}\label{psi-2-bounds}
(\varphi(x,t),\dot\varphi(x,t))\in C\sb{b}(\R,\E).
\end{equation}
The Fourier transform
\begin{equation}
\hat\varphi(x,\omega)
=\mathcal{F}\sb{t\to\omega}[\varphi(x,t)],
\qquad
(x,\omega)\in\R^2,
\end{equation}
is a continuous function of $x\in\R$ 
with 
values in 
tempered distributions of $\omega\in\R$.
It satisfies the following equation (Cf. (\ref{kg-no-cp-2})):
\begin{equation}\label{kg-no-cp-2-st}
-\omega^2\hat\varphi(x,\omega)
=
\hat\varphi''(x,\omega)
-m^2\hat\varphi(x,\omega)
+\sum\sb{J}\delta(x-X\sb{J})\hat f\sb{J}(\omega),
\qquad
(x,\omega)\in\R^2.
\end{equation}

We are going to construct a representation
for the solution $\hat\varphi(x,\omega)$
in a form suitable for our purposes.

\begin{lemma}\label{lemma-41}
$\hat\varphi$ is a smooth function of $x\in\R\backslash\mathcal{X}$
(where $\mathcal{X}=\{X\sb 1,\,X\sb 2,\,\dots,\,X\sb{N}\}$),
with values in tempered distributions of $\omega\in\R$,
and there exist quasimeasures
$\hat\varPhi\sb{J}\sp{\pm}$,
$1\le J\le N$,
and
$\hat\varTheta\sb J$,
$1\le J\le N-1$,
so that
\begin{equation}\label{c-s}
\hat\varphi(x,\omega)
=
\left\{\begin{array}{l}
\hat\varPhi\sb{1}\sp{+}(\omega)e^{-ik\sb{+}(\omega)(x-X\sb{1})}
+
\hat\varPhi\sb{1}\sp{-}(\omega)e^{-ik\sb{-}(\omega)(x-X\sb{1})},
\qquad
x\le X\sb{1},
\\
\hat\varPhi\sb{J}(\omega)
\cos(k\sb{+}(\omega)(x-X\sb{J}))
+\hat\varTheta\sb{J}(\omega)\frac{\sin(k\sb{+}(\omega)(x-X\sb{J}))}{k\sb{+}(\omega)},
\quad
x\in [X\sb{J},X\sb{J+1}],\quad 1\le J\le N-1,
\\
\hat\varPhi\sb{N}\sp{+}(\omega)e^{ik\sb{+}(\omega)(x-X\sb{N})}
+
\hat\varPhi\sb{N}\sp{-}(\omega)e^{ik\sb{-}(\omega)(x-X\sb{N})},
\qquad
x\ge X\sb{N},
\end{array}
\right.
\end{equation}
where
$
\hat\varPhi\sb{J}(\omega)
:=
\hat\varPhi\sb{J}\sp{+}(\omega)+\hat\varPhi\sb{J}\sp{-}(\omega).
$
\end{lemma}

\begin{remark}\label{remark-quasimeasure}
A tempered distribution $\mu(\omega)\in\mathscr{S}'(\R)$
is called a {\it quasimeasure} if
$\check\mu(t)=\mathscr{F}\sp{-1}\sb{\omega\to t}[\mu(\omega)]\in C\sb{b}(\R)$.
For more details, see \cite[Appendix B]{ubk-arma}.
\end{remark}

\begin{remark}\label{remark-relations}
The representation (\ref{c-s})
implies that
\begin{equation}\label{a-is-varphi}
\hat\varPhi\sb{J}(\omega)
=
\hat\varphi(X\sb J,\omega),
\quad
1\le J\le N,
\end{equation}
\begin{equation}\label{varphi-is-varphi}
\hat\varPhi\sb{1}\sp{+}(\omega)+\hat\varPhi\sb{1}\sp{-}(\omega)
=\hat\varPhi\sb{1}(\omega)=\hat\varphi(X\sb{1},\omega),
\qquad
\hat\varPhi\sb{N}\sp{+}(\omega)+\hat\varPhi\sb{N}\sp{-}(\omega)
=\hat\varphi(X\sb{N},\omega),
\end{equation}
and also that
\begin{equation}\label{a-is-varphi-prime}
\hat\varphi'(X\sb J+0,\omega)=\hat\varTheta\sb{J}(\omega),
\quad
1\le J\le N-1.
\end{equation}
\end{remark}

\begin{proof}
{\bf Step 1: Complex Fourier-Laplace transform.}\quad
We denote
\begin{equation}\label{def-f-pm}
f\sp{\pm}\sb{J}(t):=\theta(\pm t)f\sb{J}(t)
=\theta(t)F\sb{J}(\psi(X\sb{J},t))
\end{equation}
and split $\varphi(x,t)$
into
\begin{equation}\label{def-varphi-pm}
\varphi(x,t)=\varphi\sp{+}(x,t)+\varphi\sp{-}(x,t),
\qquad
{\rm where}
\quad
\varphi\sp\pm(x,t):=\theta(\pm t)\varphi(x,t).
\end{equation}
Then $\varphi\sp\pm(x,t)$ satisfy
\begin{equation}
\ddot\varphi\sp\pm(x,t)=\p\sb x^2\varphi\sp\pm(x,t)-m^2\varphi\sp\pm(x,t)
+
\sum\sb{J}
\delta(x-X\sb{J})f\sb{J}\sp\pm(t),
\qquad
t\in\R,
\label{kg-no-cp-2-pm}
\end{equation}
since $(\varphi\sp\pm,\dot\varphi\sp\pm)\at{t=0}=(0,0)$.
Let us analyze the complex Fourier-Laplace transforms of
$\varphi\sp\pm(x,t)$:
\begin{equation}\label{FL}
\displaystyle
\tilde\varphi\sp\pm(x,\omega)
=\mathcal{F}\sb{t\to\omega}[\theta(\pm t)\varphi(x,t)]
:=
\int\sb{0}\sp\infty e\sp{i\omega t}\theta(\pm t)\varphi(x,t)\,dt,
\quad\omega\in\C\sp{\pm},
\end{equation}
where
$\C\sp{\pm}:=\{z\in\C:\;\pm\Im z>0\}$.
Due to (\ref{psi-2-bounds}),
$\tilde\varphi\sp{\pm}(\cdot,\omega)$
are $H\sp{1}$-valued analytic functions of $\omega\in\C\sp{\pm}$.
In what follows, we will consider $\varphi\sp{+}$;
the function $\varphi\sp{-}$ considered in the same way.

Equation (\ref{kg-no-cp-2-pm})
implies that
$\tilde\varphi\sp{+}$ satisfies
\begin{equation}\label{FL-2}
-\omega\sp 2\tilde\varphi\sp{+}(x,\omega)
=
\p\sb x^2\tilde\varphi\sp{+}(x,\omega)-m^2\tilde\varphi\sp{+}(x,\omega)
+\sum\sb{J}\delta(x-X\sb{J})\tilde f\sb{J}\sp{+}(\omega),
\quad\omega\in\C\sp{+}.
\end{equation}
The fundamental solutions
$\displaystyle
G\sb\pm(x,\omega)=\frac{e\sp{\pm i k(\omega)\abs{x}}}{\pm 2i k(\omega)}$
satisfy
\[
G\sb\pm''(x,\omega)+(\omega^2-m^2)G\sb\pm(x,\omega)
=\delta(x),
\qquad\omega\in\C\sp{+}.
\]
The solution $\tilde\varphi\sp{+}(x,\omega)$
could be written as a linear combination of these fundamental solutions.
We use the standard ``limiting absorption principle''
for the selection of the appropriate fundamental solution: Since
$\tilde\varphi\sp{+}(\cdot,\omega)\in H\sp{1}$
for $\omega\in\C\sp{+}$,
only $G\sb{+}$ is acceptable,
because for $\omega\in\C\sp{+}$ the function $G\sb{+}(\cdot,\omega)$
is in $H\sp{1}$
by definition (\ref{def-k}),
while $G\sb{-}$ is not.
This suggests the following representation:
\begin{equation}\label{tilde-psi-tilde-f}
\tilde\varphi\sp{+}(x,\omega)
=-\sum\sb{J}
\tilde f\sp{+}\sb{J}(\omega)G\sb{+}(x-X\sb{J},\omega)
=-\sum\sb{J}
\tilde f\sp{+}\sb{J}(\omega)\frac{e\sp{i k(\omega)\abs{x-X\sb{J}}}}{2i k(\omega)},
\qquad\omega\in\C\sp{+}.
\end{equation}
The proof is straightforward since (\ref{tilde-psi-tilde-f})
belongs to $H\sp{1}(\R)$ for $\omega\in\C\sp{+}$
while the solution to (\ref{FL-2})
which is an $H\sp{1}$-valued analytic function in $\omega$
is unique.
For $x\le X\sb{1}$, the relation (\ref{tilde-psi-tilde-f}) yields
\begin{equation}\label{phi-1}
\tilde\varphi\sp{+}(x,\omega)
=-\sum\sb{J}
\tilde f\sp{+}\sb{J}(\omega)\frac{e\sp{-i k(\omega)(x-X\sb{J})}}{2i k(\omega)}
=e^{-i k(\omega)(x-X\sb{1})}\tilde\varphi\sp{+}(X\sb{1},\omega),
\qquad
x\le X\sb{1},
\quad
\omega\in\C\sp{+}.
\end{equation}
For $x\in[X\sb{J},X\sb{J+1}]$, $1\le J\le N-1$,
the relation (\ref{tilde-psi-tilde-f}) implies that
\begin{equation}\label{ast0}
\tilde\varphi\sp{+}(x,\omega)
=\tilde\varPhi\sb{J}\sp{+}(\omega)\cos(k(\omega)(x-X\sb{J}))
+\tilde\varTheta\sb{J}\sp{+}(\omega)
\frac{\sin(k(\omega)(x-X\sb{J}))}{k(\omega)},
\qquad
x\in[X\sb{J},X\sb{J+1}],
\quad
\omega\in\C\sp{+},
\end{equation}
where $\tilde\varPhi\sb{J}\sp{+}$ and $\tilde\varTheta\sb{J}\sp{+}$,
$1\le J\le N-1$,
are analytic functions of $\omega\in\C\sp{+}$.
We note that,
by (\ref{tilde-psi-tilde-f}),
\begin{equation}\label{a-b}
\tilde\varPhi\sb{J}\sp{+}(\omega)
=\tilde\varphi\sp{+}(X\sb{J},\omega),
\qquad
\tilde\varTheta\sb{J}\sp{+}(\omega)
=\p\sb{x}\tilde\varphi\sp{+}(X\sb{J}+0,\omega)
=-\sum\sb{J'}
\sgn(X\sb{J}-X\sb{J'})
\tilde f\sb{J'}\sp{+}(\omega)\frac{e^{ik(\omega)\abs{X\sb{J}-X\sb{J'}}}}{2}.
\end{equation}

\noindent
{\bf Step 2: Traces on real line.}\quad
Now we need to extend the relations (\ref{phi-1}) and (\ref{ast0})
to $\omega\in\R$.
The Fourier transform
$\hat\varphi\sp{+}(x,\omega)
:=\mathcal{F}\sb{t\to\omega}[\theta(t)\varphi(x,t)]$
is a tempered $H\sp{1}$-valued distribution of $\omega\in\R$
by (\ref{psi-2-bounds}).
It is the boundary value
of the analytic function $\tilde\varphi\sp{+}(x,\omega)$,
in the following sense:
\begin{equation}\label{bvp1}
\hat\varphi\sp{+}(x,\omega)
=\lim\limits\sb{\varepsilon\to 0+}\tilde\varphi\sp{+}(x,\omega+i\varepsilon),
\qquad\omega\in\R,
\end{equation}
where the convergence
is in the space
of tempered distributions
$\mathscr{S}'(\R,H\sp{1}(\R))$.
Indeed,
\[
\tilde\varphi\sp{+}(x,\omega+i\varepsilon)
=\mathcal{F}\sb{t\to\omega}[\theta(t)\varphi(x,t)e\sp{-\varepsilon t}],
\qquad
\theta(t)\varphi(x,t)e\sp{-\varepsilon t}
\mathop{\longrightarrow}\limits\sb{\varepsilon\to 0+}
\theta(t)\varphi(x,t),
\]
where the convergence holds in $\mathscr{S}'(\R,H\sp{1}(\R))$.
Therefore, (\ref{bvp1}) holds by the continuity of the Fourier transform
$\mathcal{F}\sb{t\to\omega}$ in $\mathscr{S}'(\R)$.

The distributions
$\hat\varPhi\sb{J}\sp{+}(\omega)$,
$\hat\varTheta\sb{J}\sp{+}(\omega)\in\mathscr{S}'(\R)$,
$\omega\in\R$,
are
defined as
the boundary values of the
functions
$\tilde\varPhi\sb{J}\sp{+}(\omega)$ and $\tilde\varTheta\sb{J}\sp{+}(\omega)$
analytic in $\omega\in\C\sp{+}$:
\begin{eqnarray}\label{phi-to-phi}
&&
\hat\varPhi\sb{J}\sp{+}(\omega)=\lim\limits\sb{\varepsilon\to 0+}
\tilde\varPhi\sb{J}\sp{+}(\omega+i\varepsilon),
\qquad
\omega\in\R,
\quad
0\le J\le N,
\\
&&
\label{theta-to-theta}
\hat\varTheta\sb{J}\sp{+}(\omega)
=\lim\limits\sb{\varepsilon\to 0+}
\tilde\varTheta\sb{J}\sp{+}(\omega+i\varepsilon),
\qquad
\omega\in\R,
\quad
1\le J\le N-1.
\end{eqnarray}
The above convergence holds in the space of quasimeasures
by (\ref{a-b}),
since
$\tilde\varphi\sp{+}(X\sb{J},\omega)$
and $\tilde f\sb{J}\sp{+}(\omega)$
are quasimeasures (see Remark~\ref{remark-quasimeasure})
while the exponential factors
in (\ref{a-b})
are multiplicators
in the space of quasimeasures \cite[Appendix B]{ubk-arma}.
Therefore, the
formulas (\ref{ast0}) with $1\le J\le N-1$
imply, in the limit $\Im\omega\to 0+$, that
\begin{equation}\label{ast1}
\hat\varphi\sp{+}(x,\omega)
=\hat\varPhi\sb{J}\sp{+}(\omega)\cos(k(\omega+i0)(x-X\sb{J}))
+\hat\varTheta\sb{J}\sp{+}(\omega)\frac{\sin(k(\omega+i0)(x-X\sb{J}))}{k(\omega+i0)},
\qquad
x\in [X\sb{J},X\sb{J+1}],
\quad
\omega\in\R,
\end{equation}
since
$\cos(k(\omega+i0)(x-X\sb{J}))$
and
$\frac{\sin(k(\omega+i0)(x-X\sb{J}))}{k(\omega+i0)}$
are smooth functions of $\omega\in\R$.
Similar representation holds for $\hat\varphi\sp{-}(x,\omega)$.
Therefore, the representation (\ref{c-s}) follows
for $X\sb{1}\le x\le X\sb{N}$.

The formula (\ref{c-s}) for $x\le X\sb{1}$
follows from taking the limit
$\Im\omega\to 0+$
in the expression (\ref{phi-1})
for $\tilde\varphi\sp{+}(x,\omega)$
and the limit
$\Im\omega\to 0-$
in a similar expression
for $\tilde\varphi\sp{-}(x,\omega)$:
\begin{equation}\label{phi-1-m}
\tilde\varphi\sp{-}(x,\omega)
=-\sum\sb{J}
\tilde f\sp{-}\sb{J}(\omega)\frac{e\sp{-i k(\omega)(x-X\sb{J})}}{2i k(\omega)}
=e^{-i k(\omega)(x-X\sb{1})}\tilde\varphi\sp{-}(X\sb{1},\omega),
\qquad
x\le X\sb{1},
\quad
\omega\in\C\sp{-},
\end{equation}
and then taking the sum of the resulting expressions.
This justifies (\ref{c-s}) for $x\le X\sb{1}$.
Similarly we justify (\ref{c-s}) for $x\ge X\sb{N}$.
\end{proof}

\section{Absolute continuity of the spectrum}
\label{sect-ac}

\begin{lemma}\label{lemma-continuity-0}
The distributions
$\hat\varPhi\sb{1}\sp\pm(\omega)$, $\hat\varPhi\sb{N}\sp\pm(\omega)$
are absolutely continuous for $\abs{\omega}>m$,
and moreover
\begin{equation}\label{phi-0-41}
\int\sb{\abs{\omega}>m}
\left[
\abs{\hat\varPhi\sb{1}\sp\pm(\omega)}\sp 2
+
\abs{\hat\varPhi\sb{N}\sp\pm(\omega)}\sp 2
\right]
\omega\,k\sb{+}(\omega)
\,d\omega<\infty,
\end{equation}
where $\omega\,k\sb{+}(\omega)\ge 0$ by {(\ref{omega-k})}.
\end{lemma}

The bound for each of
$\hat\varPhi\sb{1}\sp\pm(\omega)$,
$\hat\varPhi\sb{N}\sp\pm(\omega)$
is obtained verbatim by applying the proof of
\cite[Proposition 3.3]{ubk-arma}.

\begin{proposition}\label{prop-continuity}
The distributions
$\hat\varPhi\sb{J}(\omega)$, $1\le J\le N$,
and
$\hat\varTheta\sb{J}(\omega)$, $1\le J\le N-1$,
are absolutely continuous for
$\abs{\omega}>\mu\sb{J}$ and $\abs{\omega}>(2p\sb{J}-1)\mu\sb{J}$,
respectively,
with $\mu\sb{J}$ defined in (\ref{def-mu}).
Moreover,
for any $\epsilon>0$,
\begin{equation}\label{41}
\int\sb{\abs{\omega}
>\mu\sb{J}+\epsilon}
\abs{\hat\varPhi\sb{J}(\omega)}\sp 2
\omega^2
\,d\omega<\infty,
\quad
1\le J\le N;
\qquad
\int\sb{\abs{\omega}>(2p\sb{J}-1)\mu\sb{J}+\epsilon}
\abs{\hat\varTheta\sb{J}(\omega)}\sp 2
\,d\omega<\infty,
\quad
1\le J\le N-1.
\end{equation}
\end{proposition}

\begin{proof}
We will use induction,
proving the absolute continuity of
$\hat\varphi(X\sb{J},\omega)$
and $\p\sb x\hat\varphi(X\sb{J}\pm 0,\omega)$
starting with $J=1$ and going to $J=N$.
By Lemma~\ref{lemma-41},
$\hat\varphi(X\sb 1,\omega)
=\hat\varPhi\sb{1}(\omega)
=\hat\varPhi\sb{1}\sp{+}(\omega)+\hat\varPhi\sb{1}\sp{-}(\omega)$
and
$
\p\sb x\hat\varphi(X\sb 1-0,\omega)=
-i k\sb{+}(\omega)\hat\varPhi\sb{1}\sp{+}(\omega)
-i k\sb{-}(\omega)\hat\varPhi\sb{1}\sp{-}(\omega).
$
Hence,
Lemma~\ref{lemma-continuity-0}
implies that,
for any $\epsilon>0$,
\begin{equation}\label{phi-phi-1}
\int\sb{\abs{\omega}>m+\epsilon}
\abs{\hat\varphi(X\sb{1},\omega)}\sp 2
\omega^2
\,d\omega<\infty,
\qquad
\int\sb{\abs{\omega}>m+\epsilon}
\abs{\hat\varphi'(X\sb{1}-0,\omega)}\sp 2
\,d\omega<\infty.
\end{equation}
Now assume that for some $1\le J<N$ and for any $\epsilon>0$
we have:
\begin{equation}\label{phi-phi-j}
\int\sb{\abs{\omega}>\mu\sb{J}+\epsilon}
\abs{\hat\varphi(X\sb{J},\omega)}\sp 2
\omega^2
\,d\omega<\infty,
\qquad
\int\sb{\abs{\omega}>\mu\sb{J}+\epsilon}
\abs{\hat\varphi'(X\sb{J}-0,\omega)}\sp 2
\,d\omega<\infty.
\end{equation}
Lemma~\ref{lemma-41}
and equation (\ref{kg-no-cp-2-st})
yield the jump condition
\begin{equation}\label{jump-j}
\hat\varTheta\sb{J}(\omega)
=
\hat\varphi'(X\sb{J}+0,\omega)
=\hat\varphi'(X\sb{J}-0,\omega)
-\hat f\sb{J}(\omega),
\qquad\omega\in\R,
\end{equation}
where $f\sb{J}(t)=F\sb{J}(\psi(X\sb{J},t))$ by (\ref{def-f}).

\begin{lemma}\label{lemma-fj-ac}
For any $\epsilon>0$
the following inequality holds:
\begin{equation}\label{int-f2}
\int\sb{\abs{\omega}>(2p\sb{J}-1)(\mu\sb{J}+2\epsilon)}
\abs{\hat f\sb{J}(\omega)}^2
\,d\omega<\infty.
\end{equation}
\end{lemma}

\begin{proof}
Let
$\zeta\sb{J}(\omega)\in C\sp\infty\sb{0}(\R)$
be such that
$\zeta\sb{J}(\omega)\equiv 1$ for $\abs{\omega}\le\mu\sb{J}+\epsilon$
and
$\zeta\sb{J}(\omega) \equiv 0$ for $\abs{\omega}\ge\mu\sb{J}+2\epsilon$.
We denote $\psi(X\sb{J},t)$ by $\uppsi\sb{J}(t)$,
and split it into
\begin{equation}\label{psi-psi-psi}
\uppsi\sb{J}(t)
=\uppsi\sb{J,b}(t)+\uppsi\sb{J,d}(t),
\end{equation}
where the functions in the right-hand side are defined
by their Fourier transforms:
\begin{equation}
\hat\uppsi\sb{J,b}(\omega)
=\zeta\sb{J}(\omega)\hat\uppsi\sb{J}(\omega)
=\zeta\sb{J}(\omega)\hat\psi(X\sb{J},\omega),
\qquad
\hat\uppsi\sb{J,d}(\omega)
=(1-\zeta\sb{J}(\omega))\hat\uppsi\sb{J}(\omega)
=(1-\zeta\sb{J}(\omega))\hat\psi(X\sb{J},\omega).
\end{equation}
By Lemma~\ref{lemma-chi-h1}
and by (\ref{phi-phi-j}),
we have
\begin{equation}\label{chi-varphi-h1}
\int\sb{\R}
\Abs{(1-\zeta\sb{J}(\omega))
\hat\chi(X\sb{J},\omega)}^2\omega^2\,d\omega<\infty,
\qquad
\int\sb{\R}
\Abs{(1-\zeta\sb{J}(\omega))
\hat\varphi(X\sb{J},\omega)}^2\omega^2\,d\omega
<\infty.
\end{equation}
Since
$
\hat\uppsi\sb{J,d}(\omega)
=(1-\zeta\sb{J}(\omega))
(\hat\chi(X\sb{J},\omega)+\hat\varphi(X\sb{J},\omega)),
$
we also have
\[
\int\sb{\R}
\Abs{(1-\zeta\sb{J}(\omega))
\hat\uppsi\sb{J}(\omega)}^2\omega^2\,d\omega
<\infty,
\]
proving that
\begin{equation}\label{psi-jd-h1}
\uppsi\sb{J,d}(t)\in H\sp 1(\R).
\end{equation}
For
$
\hat f\sb{J}(\omega)
=
\mathcal{F}\sb{t\to\omega}[F\sb{J}(\uppsi\sb{J}(t))]
=
\mathcal{F}\sb{t\to\omega}[F\sb{J}(\psi(X\sb{J},t))]
$,
taking into account (\ref{f-is-such}) and (\ref{psi-psi-psi}),
we have:
\begin{eqnarray}
\hat f\sb{J}(\omega)
&=&-\sum\sb{n=1}\sp{p\sb{J}}
2n\,u\sb{J,n}
\underbrace
{
(\hat\uppsi\sb{J}\ast\hat{\overline\uppsi}\sb{J})
\ast\dots\ast
(\hat\uppsi\sb{J}\ast\hat{\overline\uppsi}\sb{J})
}\sb{n-1}
\ast
\hat\uppsi\sb{J}
\nonumber
\\
&=&
.\ .\ .\ .\ .\ \ -\sum\sb{n=1}\sp{p\sb{J}}
2n\,u\sb{J,n}
\underbrace
{
(\hat\uppsi\sb{J,b}\ast\hat{\overline\uppsi}\sb{J,b})
\ast\dots\ast
(\hat\uppsi\sb{J,b}\ast\hat{\overline\uppsi}\sb{J,b})
}\sb{n-1}
\ast
\hat\uppsi\sb{J,b},
\label{conv-conv}
\end{eqnarray}
where the dots in the right-hand side
denote the convolutions of
$\hat\uppsi\sb{J,b}$, $\hat{\overline\uppsi}\sb{J,b}$,
$\hat\uppsi\sb{J,d}$, and $\hat{\overline\uppsi}\sb{J,d}$
that contain at least one of
$\hat\uppsi\sb{J,d}$, $\hat{\overline\uppsi}\sb{J,d}$.
Since
$\uppsi\sb{J,b}(t)$, $\uppsi\sb{J,d}(t)$ are bounded
while $\uppsi\sb{J,d}(t)\in H\sp 1(\R)$
by (\ref{psi-jd-h1}),
all these terms belong to $L\sp 2(\R)$.
Finally, since
$\supp\hat\uppsi\sb{J,b}
\subset[-\mu\sb{J}-2\epsilon,\mu\sb{J}+2\epsilon]$,
the convolutions
under the summation sign
in the right-hand side of (\ref{conv-conv})
are supported inside
$[-(2p\sb{J}-1)(\mu\sb{J}+2\epsilon),(2p\sb{J}-1)(\mu\sb{J}+2\epsilon)]$
and do not contribute into the integral (\ref{int-f2}).
\end{proof}

Using
(\ref{phi-phi-j}) and Lemma~\ref{lemma-fj-ac}
to estimate the norms of
$\p\sb x\hat\varphi(X\sb{J}-0,\omega)$ and $\hat f\sb{J}(\omega)$
in the right-hand side in the relation (\ref{jump-j}),
we conclude that
\begin{equation}\label{bj-ac}
\int\limits\sb{\abs{\omega}>(2p\sb{J}-1)(\mu\sb{J}+2\epsilon)}
\Abs{
\hat\varphi'(X\sb{J}+0,\omega)}^2\,d\omega
<\infty.
\end{equation}
Now the inequalities
\begin{equation}\label{phi-phi-jj}
\int\limits\sb{\abs{\omega}>(2p\sb{J}-1)(\mu\sb{J}+2\epsilon)}
\abs{\hat\varphi(X\sb{J+1},\omega)}^2
\omega^2\,d\omega<\infty,
\qquad
\int\limits\sb{\abs{\omega}>(2p\sb{J}-1)(\mu\sb{J}+2\epsilon)}
\abs{\hat\varphi'(X\sb{J+1}-0,\omega)}^2
\,d\omega<\infty
\end{equation}
follow
from the representation (\ref{c-s}) for $x\in[X\sb{J},X\sb{J+1}]$,
where we apply the first inequality from (\ref{phi-phi-j})
and the inequality (\ref{bj-ac}).
Therefore, starting with (\ref{phi-phi-1}),
one shows by induction that (\ref{phi-phi-j}) holds for all $1\le J\le N$.
The estimates on
$\hat\varPhi\sb{J}(\omega)=\hat\varphi(X\sb{J},\omega)$
and 
$\hat\varTheta\sb{J}(\omega)=\hat\varphi'(X\sb{J}+0,\omega)$
stated in the Proposition
follow from (\ref{phi-phi-j}) and (\ref{bj-ac}), respectively.
This finishes the proof of Proposition~\ref{prop-continuity}.
\end{proof}

\begin{corollary}\label{corollary-continuity}
The distributions
$\hat\varPhi\sb{J}(\omega)=\hat\varphi(X\sb J,\omega)$,
$1\le J\le N$,
are absolutely continuous for
$\abs{\omega}>{M}\sb{J}$,
while $\hat\varTheta\sb{J}(\omega)=\p\sb x\hat\varphi(X\sb J+0,\omega)$,
$1\le J\le N-1$,
are absolutely continuous for $\abs{\omega}>(2p\sb{J}-1){M}\sb{J}$,
where ${M}\sb{J}:=\min(\mu\sb{J},\mu'\sb{J})$
is defined in (\ref{def-m-lambda}).
\end{corollary}

\begin{proof}
In the proof of Proposition~\ref{prop-continuity},
we could as well proceed from $J=N$ to $J=1$,
proving the result stated in the Corollary.
\end{proof}

\section{Compactness}
\label{sect-bound}

\subsection*{Second dispersive component}
Let $\zeta(\omega)\in C\sp\infty\sb{0}(\R)$ be such that
$\zeta(\omega)\equiv 1$ for $\abs{\omega}<\Lambda$,
where $\Lambda$ is from (\ref{def-m-lambda}).
Define $\varphi\sb{d}(x,t)$
by its Fourier transform:
\begin{equation}\label{dd}
\hat\varphi\sb{d}(x,\omega)
:=
(1-\zeta(\omega))\hat\varphi(x,\omega)
\ \ x\in\R,\ \ \omega\in\R.
\end{equation}

\begin{lemma}\label{lemma-phi-d-b}
$\varphi\sb{d}(x,t)$ is a bounded continuous
function
of $t\in\R$
with values in $H\sp 1(\R)$:
\begin{equation}\label{phi-d-b}
\varphi\sb{d}(x,t)\in C\sb{b}(\R,H\sp{1}(\R)).
\end{equation}
The local energy decay holds for $\varphi\sb{d}(x,t)$:
\begin{equation}\label{led}
\lim\sb{t\to\infty}
\norm{(\varphi\sb{d},\dot\varphi\sb{d})}\sb{\E,R}=0,
\qquad\forall R>0.
\end{equation}
\end{lemma}

\begin{proof}
We generalize the proof of \cite[Proposition 3.6]{ubk-arma}.
By Lemma~\ref{lemma-41},
\begin{equation}\label{ftb00}
\hat\varphi\sb{d}(x,\omega)
=
\left\{
\begin{array}{l}
(1-\zeta(\omega))
\left[
\hat\varPhi\sb{1}\sp{+}(\omega)e\sp{-ik\sb{+}(\omega)(x-X\sb{1})}
+
\hat\varPhi\sb{1}\sp{-}(\omega)e\sp{-ik\sb{-}(\omega)(x-X\sb{1})}
\right]
,
\qquad
x\le X\sb{1},
\\
(1-\zeta(\omega))\hat\varPhi\sb{J}(\omega)\cos(k\sb{+}(\omega)(x-X\sb{J}))
+
(1-\zeta(\omega))
\hat\varTheta\sb{J}(\omega)
\frac{\sin(k\sb{+}(\omega)(x-X\sb{J}))}{k\sb{+}(\omega)(x-X\sb{J})},
\qquad
x\in[X\sb{J},X\sb{J+1}],
\\
(1-\zeta(\omega))
\left[
\hat\varPhi\sb{N}\sp{+}(\omega)e\sp{ik\sb{+}(\omega)(x-X\sb{N})}
+
\hat\varPhi\sb{N}\sp{-}(\omega)e\sp{ik\sb{-}(\omega)(x-X\sb{N})}
\right],
\qquad
x\ge X\sb{N}.
\end{array}
\right.
\end{equation}
Each of the functions entering the above expression,
considered on the whole real line,
corresponds to a finite energy solution to a linear
Klein-Gordon equation,
satisfying the properties stated in the lemma.
For example,
define $u(x,t)$ by its Fourier transform:
\[
\hat u(x,\omega)
:=
(1-\zeta(\omega))\hat\varPhi\sb{1}(\omega)
\cos(k\sb{+}(\omega)(x-X\sb{1})),
\qquad
x\in\R.
\]
Then $u(x,t)$ is a solution to a linear Klein-Gordon equation,
and, by Proposition~\ref{prop-continuity},
the corresponding initial data are of finite energy:
\[
(u(x,0),\dot u(x,0))\in\E.
\]
Hence
$u(x,t)\in C\sb{b}(\R,H\sp 1(\R))$
and satisfies the local energy decay of the form (\ref{led})
(see \cite[Lemma 3.1]{ubk-arma}.
This finishes the proof.
\end{proof}

\subsection*{Compactness for the bound component}

We introduce the bound component of $\varphi(x,t)$
by
\begin{equation}\label{bb}
\varphi\sb{b}(x,t)
=\varphi(x,t)-\varphi\sb{d}(x,t)
=\psi(x,t)-\chi(x,t)-\varphi\sb{d}(x,t),
\qquad
x\in\R,\ \ t\in\R.
\end{equation}
By Lemma~\ref{lemma-phi-d-b},
\begin{equation}\label{ebb}
\varphi\sb{b}(x,t)\in C\sb{b}(\R,H\sp{1}(\R)).
\end{equation}
Lemma~\ref{lemma-41} and (\ref{dd}), (\ref{bb})
imply the multiplicative relation
\begin{equation}\label{ftb0}
\hat\varphi\sb{b}(x,\omega)
=
\left\{
\begin{array}{l}
\zeta(\omega)
\big[
\hat\varPhi\sb{1}\sp{+}(\omega)e\sp{-i k\sb{+}(\omega)(x-X\sb{1})}
+
\hat\varPhi\sb{1}\sp{-}(\omega)e\sp{-i k\sb{-}(\omega)(x-X\sb{1})}
\big],
\qquad
x\le X\sb{1},
\\
\zeta(\omega)
\big[\hat\varPhi\sb{J}(\omega)\cos(k\sb{+}(\omega)(x-X\sb{J}))
+
\hat\varTheta\sb{J}(\omega)
\frac{\sin(k\sb{+}(\omega)(x-X\sb{J}))}{k\sb{+}(\omega)}
\big],
\qquad
x\in[X\sb{J},X\sb{J+1}],
\\
\zeta(\omega)
\big[
\hat\varPhi\sb{N}\sp{+}(\omega)e\sp{i k\sb{+}(\omega)(x-X\sb{N})}
+
\hat\varPhi\sb{N}\sp{-}(\omega)e\sp{i k\sb{-}(\omega)(x-X\sb{N})}
\big],
\qquad
x\ge X\sb{N}.
\end{array}
\right.
\end{equation}
By (\ref{ebb}), the functions
\[
\varphi\sb{b,J}(t)
:=
\varphi\sb{b}(X\sb{J},t)
=\varphi(X\sb{J},t)-\varphi\sb{d}(X\sb{J},t)
\]
are bounded and continuous.
Therefore,
$
\hat\varphi\sb{b}(X\sb{J},\cdot)\in\mathscr{S}'(\R)
$
are quasimeasures (see Remark~\ref{remark-quasimeasure}).

\begin{proposition}
\begin{enumerate}
\item
The function $\varphi\sb{b}(x,t)$ is smooth
for $x\in\R\backslash\mathcal{X}$
(where $\mathcal{X}=\{X\sb 1,\,X\sb 2,\,\dots,\,X\sb{N}\}$)
and $t\in\R$.
\item
For any $R>0$,
\begin{equation}\label{bqda}
\sup\limits\sb{
\abs{x}\le R,\,x\notin\mathcal{X}}
\,\,
\sup\limits\sb{t\in\R}
\abs{\p\sb x\sp\xx\p\sb t\sp\yy\varphi\sb{b}(x,t)}
<\infty.
\end{equation}
\end{enumerate}
\end{proposition}

The argument repeats
the proof of Proposition~\cite[Proposition 4.1]{ubk-arma}.

\begin{remark}
Let us note that
the bounds (\ref{bqda}) are independent of $x$
and remain valid for
$x\notin\mathcal{X}$,
although the derivatives
$\p\sb x\sp\xx\p\sb t\sp\yy\varphi\sb{b}(x,t)$
with $\xx\ne 0$
may have jumps at $x=X\sb{J}$.
(Note that
this is the case for the solitary waves in (\ref{solitary-wave-profile}).)
\end{remark}

We now may deduce the compactness of the set of
translations of the bound component,
$\{\varphi\sb{b}(x,s+t)\sothat s\ge 0\}$.

\begin{corollary}\label{coco}
\begin{enumerate}
\item
By the Ascoli-Arzel\`a Theorem,
for any sequence $s\sb{j}\to\infty$
there exists a subsequence $s\sb{j'}\to\infty$
such that
\begin{equation}
\varphi\sb{b}(x,s\sb{j'}+t)
\to
\beta(x,t),
\qquad
x\in\R,
\quad
t\in\R,
\label{ol}
\end{equation}
and also for any nonnegative integers $\xx$ and $\yy$,
\begin{equation}
\p\sb x\sp\xx\p\sb t\sp\yy
\varphi\sb{b}(x,s\sb{j'}+t)
\to
\p\sb x\sp\xx\p\sb t\sp\yy
\beta(x,t),
\qquad
x\notin\mathcal{X},
\quad
t\in\R,
\label{olpd}
\end{equation}
for some
$\beta(x,t)\in C\sb{b}(\R,H\sp{1}(\R))$.
The convergence in (\ref{ol}) and (\ref{olpd})
is uniform in $x$ and $t$
as long as $\abs{x}+\abs{t}\le R$, for any $R>0$.
The convergence in (\ref{olpd})
also holds for $x=X\sb{J}\pm 0$.
\item
By the Fatou Lemma,
\begin{equation}\label{beta-beta}
\sup\limits\sb{t\in\R}\norm{\beta(\cdot,t)}\sb{H\sp{1}}<\infty.
\end{equation}
\end{enumerate}
\end{corollary}

We call {\it omega-limit trajectory}
any function $\beta(x,t)$
that can appear as a limit in (\ref{ol}), (\ref{olpd}).

\begin{remark}\label{remark-end}
Previous analysis demonstrates that the long-time asymptotics
of the solution $\psi(x,t)$ in $\E\sb{F}$
depends only on the singular component
$\varphi(x,t)$.
Due to Corollary~\ref{coco},
to conclude the proof of Theorem~\ref{main-theorem},
it suffices to check that every omega-limit trajectory
belongs to the set of solitary waves;
that is,
\begin{equation}\label{eidd}
\beta(x,t)
=
\phi\sb{\omega\sb{+}}(x)e\sp{-i\omega\sb{+}t}
\qquad
{\rm for\ some\ }\omega\sb{+}\in[-m,m].
\end{equation}
\end{remark}

\section{Nonlinear spectral analysis}
\label{sect-nonlinear-analysis}
\subsection*{Bounds for the spectrum}

By Lemmas~\ref{lemma-decay-psi1} and \ref{lemma-phi-d-b},
the dispersive components $\chi(\cdot,t)$
and
$\varphi\sb{d}(\cdot,t)$ converge
to zero in
${\E}\sb{F}$
as $t\to\infty$.
On the other hand,
by Corollary~\ref{coco},
the bound component
$\varphi\sb{b}(x,t+s\sb{j'})$
converges to $\beta(x,t)$
as $j'\to\infty$,
uniformly in every compact set of the plane $\R^2$.
Hence,
$\psi(x,t+s\sb{j'})=
\varphi\sb{b}(x,t+s\sb{j'})+\chi(x,t+s\sb{j'})+\varphi\sb{d}(x,t+s\sb{j'})$
also
converges to $\beta(x,t)$
uniformly in every compact set of the plane $\R^2$.
Therefore, taking the limit in equation (\ref{kg-no}),
we conclude that
the omega-limit trajectory $\beta(x,t)$ also satisfies
the same equation:
\begin{equation}\label{kg-no-beta}
\ddot\beta(x,t)
=\beta''(x,t)-m^2\beta(x,t)+\sum\sb{J}\delta(x-X\sb{J})F\sb{J}(\beta).
\end{equation}

\begin{remark}
Note that the bound component
$\varphi\sb{b}(x,t)$ itself generally does not satisfy equation (\ref{kg-no-beta}).
\end{remark}

Taking the Fourier transform of $\beta$ in time,
we see by (\ref{olpd})
that $\hat\beta(x,\omega)$
is a continuous function of $x\in\R$,
smooth for $x\in\R\backslash\mathcal{X}$,
with values in tempered distributions of $\omega\in\R$,
and that it satisfies the corresponding stationary equation
\begin{equation}\label{kg-no-beta-1}
-\omega^2\hat\beta(x,\omega)
=\hat\beta''(x,\omega)-m^2\hat\beta(x,\omega)
+\sum\sb{J}\delta(x-X\sb{J})\hat g\sb{J}(\omega),
\qquad
(x,\omega)\in\R^2,
\end{equation}
valid in the sense of tempered distributions of $(x,\omega)\in\R^2$,
where
$\hat g\sb{J}(\omega)$
are the Fourier transforms of the functions
\begin{equation}
g\sb{J}(t):=F\sb{J}(\beta(X\sb{J},t)),
\qquad
1\le J\le N.
\end{equation}
We also denote
\begin{equation}\label{def-upbeta}
\upbeta\sb{J}(t):=\beta(X\sb{J},t),
\qquad
\Sigma\sb{J}:=\supp\hat\upbeta\sb{J},
\qquad 1\le J\le N.
\end{equation}

From (\ref{ftb0}),
we know that the spectrum of $\varphi\sb{b}(x,t)$
is bounded for all $x\in\R$.
Hence,
the convergence
(\ref{olpd}) implies that
the spectrum of $\beta(x,t)$ is also bounded.
We will need more precise bounds on the size
of the spectrum of $\beta$:

\begin{lemma}\label{lemma-m-m}
\begin{enumerate}
\item
$\Sigma\sb{J}
:=\supp\hat\upbeta\sb{J}
\subset[-{M}\sb{J},{M}\sb{J}]$,
\quad
$1\le J\le N$;
\item
$\supp\hat\beta'(X\sb{J}+0,\omega)
\subset[-(2p\sb{J}-1){M}\sb{J},(2p\sb{J}-1){M}\sb{J}]$,
\quad
$1\le J\le N-1$,
with ${M}\sb{J}>0$ defined in (\ref{def-m-lambda}).
\end{enumerate}

\end{lemma}

\begin{proof}
We have the relation
\[
\varphi\sb{b}(x,s\sb{j}+t)
=
\frac{1}{2\pi}\int\sb{\R}
e^{-i\omega t}
e^{-i \omega s\sb{j}}
\hat\varphi\sb{b}(x,\omega)\,d\omega,
\qquad
x\in\R,\quad t\in\R,
\]
where the integral is understood as the pairing
of a smooth function (oscillating exponent) with a compactly supported
distribution.
Then the convergence (\ref{ol}) implies that
\begin{equation}\label{phi-to-beta}
e^{-i\omega s\sb{j'}}
\hat\varphi\sb{b}(x,\omega)
\to\hat\beta(x,\omega),
\qquad
x\in\R,
\quad s\sb{j'}\to\infty,
\end{equation}
in the sense of quasimeasures.
Since $\hat\varphi\sb{b}(X\sb J,\omega)$
is locally $L\sp 2$ for $\abs{\omega}>{M}\sb{J}$
by Corollary~\ref{corollary-continuity},
the convergence (\ref{phi-to-beta})
at $x=X\sb J$
shows that
$\hat\upbeta\sb{J}(\omega):=\hat\beta(X\sb{J},\omega)$
vanishes for $\abs{\omega}>{M}\sb{J}$.
This proves the first statement of the lemma.

The second statement is proved similarly.
Namely, the convergence (\ref{olpd}) implies that
\begin{equation}\label{phi-to-beta-pd}
e^{-i\omega s\sb{j'}}
\p\sb x\hat\varphi\sb{b}(X\sb{J}+0,\omega)
\to\p\sb x\hat\beta(X\sb{J}+0,\omega),
\qquad
s\sb{j'}\to\infty,
\end{equation}
in the sense of quasimeasures.
Since
$\hat\varphi\sb{b}'(X\sb{J}+0,\omega)$
is locally $L\sp 2$ for $\abs{\omega}>(2p\sb{J}-1){M}\sb{J}$
by Corollary~\ref{corollary-continuity},
the convergence (\ref{phi-to-beta-pd})
shows that
$\hat\beta'(X\sb{J}+0,\omega)$
vanishes for $\abs{\omega}>(2p\sb{J}-1){M}\sb{J}$.
\end{proof}

We denote
\begin{equation}\label{def-kappa}
\kappa(\omega):=-i k\sb{+}(\omega),
\qquad
\omega\in\R,
\end{equation}
where $k\sb{+}(\omega)$ was introduced in (\ref{def-k-plus}).
We then have
$\Re\kappa(\omega)\ge 0$,
and also
\[
\kappa(\omega)=\sqrt{\omega^2-m^2}>0\quad{\rm for}\quad -m<\omega<m,
\]
in accordance  with
(\ref{solitary-wave-profile}).

\begin{proposition}\label{prop-beta}
The distribution
$\hat\beta(x,\omega)$ admits the following representation:
\begin{equation}\label{beta-beta-c}
\hat\beta(x,\omega)
=
\left\{
\begin{array}{l}
\hat\upbeta\sb{1}(\omega)
e^{\kappa(\omega)(x-X\sb{1})},
\qquad
x\le X\sb{1},
\\
\hat\upbeta\sb{J}(\omega){\cosh(\kappa(\omega)(x-X\sb{J}))}+
\hat\beta'(X\sb{J}+0,\omega)
\frac{\sinh(\kappa(\omega)(x-X\sb{J}))}{\kappa(\omega)},
\qquad
x\in[X\sb{J},X\sb{J+1}],
\quad
1\le J\le N-1,
\\
\hat\upbeta\sb{N}(\omega)
e^{-\kappa(\omega)(x-X\sb{N})},
\qquad
x\ge X\sb{N}.
\end{array}
\right.
\end{equation}
\end{proposition}

\begin{proof}
By (\ref{phi-to-beta}),
the middle line in (\ref{beta-beta-c})
follows from the representation (\ref{c-s})
since the multiplicators are smooth bounded functions of $\omega\in\R$.
Taking the limit in the first line of (\ref{c-s}),
we obtain the first line in (\ref{beta-beta-c})
since
$\Sigma\sb{1}\subset[-m,m]$
by Lemma~\ref{lemma-m-m},
while
$k\sb{+}(\omega)=k\sb{-}(\omega)=i\kappa(\omega)$
for $-m\le\omega\le m$ (Cf. (\ref{k-plus-minus}), (\ref{def-kappa})).
Similarly we explain the last line in (\ref{beta-beta-c}).
\end{proof}

\subsection*{Reduction to point spectrum}

\begin{proposition}\label{prop-one-omega}
Any omega-limit trajectory $\beta(x,t)$ 
is a solitary wave:
\[
\beta(x,t)=\phi(x)e^{-i\omega\sb{+}t}
\quad{\rm with}
\quad
\omega\sb{+}\in[-m,m]
\quad
{\rm and}
\quad
\phi(x)\in H^1(\R).
\]
\end{proposition}

\begin{proof}
The proof is based on the following lemmas.

\begin{lemma}
If $\Sigma\sb{1}=\emptyset$,
then $\beta(x,t)\equiv 0$.
\end{lemma}

\begin{proof}
According to equation (\ref{kg-no-beta-1}),
the function
$\hat\beta$
satisfies the following continuity and jump conditions
at the point $X\sb{1}$:
\begin{equation}\label{jump-condition}
\hat\beta(X\sb{1}+0,\omega)=\hat\beta(X\sb{1}-0,\omega)
=\hat\upbeta\sb{1}(\omega),
\qquad
\hat\beta'(X\sb{1}+0,\omega)=\hat\beta'(X\sb{1}-0,\omega)
+
\hat g\sb{1}(\omega),
\qquad
\omega\in\R.
\end{equation}
$\Sigma\sb{1}=\emptyset$
means that $\hat\upbeta\sb{1}(\omega)\equiv 0$,
that is, $\upbeta\sb{1}(t)\equiv 0$.
Hence,
$g\sb{1}(t)\equiv F\sb{1}(\upbeta\sb{1}(t))\equiv 0$,
and $\hat g\sb{1}(\omega)\equiv 0$.
On the other hand, the first line of (\ref{beta-beta-c})
implies that
$\hat\beta(x,\omega)\equiv 0$ for $x\le X\sb{1}$,
and in particular
$\hat\beta'(X\sb{1}-0,\omega)\equiv 0$.
Therefore, the jump condition (\ref{jump-condition}) implies that
$\hat\beta'(X\sb{1}+0,\omega)\equiv 0$.
Hence,
$\hat\beta(x,\omega)\equiv 0$ for $x\in[X\sb{1},X\sb{2}]$
by the middle line of (\ref{beta-beta-c}).
By induction, 
$\hat\upbeta\sb{J}(x,\omega)\equiv 0$.
\end{proof}

Now we consider the case $\Sigma\sb{1}\not=\emptyset$.

\begin{lemma}
If $\Sigma\sb{1}\not=\emptyset$,
then $\Sigma\sb{1}=\{\omega\sb{+}\}$
for some $\omega\sb{+}\in[-m,m]$.
\end{lemma}

\begin{proof}
By Lemma~\ref{lemma-m-m},
we know that $\Sigma\sb{1}\subset[-m,m]$.
To show that $\Sigma\sb{1}$ consists of a single point,
we assume that, on the contrary,
$\inf\Sigma\sb{1}<\sup\Sigma\sb{1}$.
By (\ref{f-is-such}),
the Fourier transform
$\hat g\sb{1}(\omega)$
of $g\sb 1(t):=F\sb{1}(\beta(X\sb{1},t))$
is given by
\begin{equation}\label{conv}
\hat g\sb{1}
=
-\sum\sb{n=1}\sp{p\sb{1}}
2n\,u\sb{1,n}
\underbrace{
(\hat\upbeta\sb{1}\ast\hat{\overline\upbeta}\sb{1})
\ast\;\dots\;\ast
(\hat\upbeta\sb{1}\ast\hat{\overline\upbeta}\sb{1})
}\sb{n-1}
\ast\hat\upbeta\sb{1}.
\end{equation}
Applying the Titchmarsh Convolution Theorem
\cite{titchmarsh} (see also \cite[p.119]{MR1400006}
and \cite[Theorem 4.3.3]{MR1065136})
to the convolutions in (\ref{conv}),
we obtain the following equalities:
\begin{eqnarray}
\label{supp-g}
&&
\inf\supp\hat g\sb{1}
=\inf\supp\hat\upbeta\sb{1}+(p\sb{1}-1)
\inf\supp(\hat\upbeta\sb{1}\ast\hat{\overline\upbeta}\sb{1})
=\inf\Sigma\sb{1}+(p\sb{1}-1)(\inf\Sigma\sb{1}-\sup\Sigma\sb{1}),
\\
&&
\sup\supp\hat g\sb{1}
=\sup\supp\hat\upbeta\sb{1}+(p\sb{1}-1)
\sup\supp(\hat\upbeta\sb{1}\ast\hat{\overline\upbeta}\sb{1})
=\sup\Sigma\sb{1}+(p\sb{1}-1)(\sup\Sigma\sb{1}-\inf\Sigma\sb{1}),
\label{supp-gg}
\end{eqnarray}
where we used the relations
$
\inf\supp\hat{\overline\upbeta}\sb 1=-\sup\supp\hat\upbeta\sb 1,
$
$
\sup\supp\hat{\overline\upbeta}\sb 1=-\inf\supp\hat\upbeta\sb 1.
$
Note that the Titchmarsh theorem is applicable
since $\supp\hat\upbeta\sb{1}$ is compact
by Lemma~\ref{lemma-m-m}.
Since we assumed that $\inf\Sigma\sb{1}<\sup\Sigma\sb{1}$,
(\ref{supp-g}) and (\ref{supp-gg})
imply that
$
\inf\supp \hat g\sb{1}
<\inf\Sigma\sb{1},
$
$
\sup\supp \hat g\sb{1}
>\sup\Sigma\sb{1}.
$
Therefore, the jump condition (\ref{jump-condition})
with $J=1$
implies that
\begin{equation}\label{imp}
\inf\supp \hat\beta'(X\sb{1}+0,\cdot)
=
\inf\supp\hat g\sb{1}
<\inf\Sigma\sb{1},
\qquad
\sup\supp \hat\beta'(X\sb{1}+0,\cdot)=
\sup\supp\hat g\sb{1}
>
\sup\Sigma\sb{1}.
\end{equation}
The ratio
$
\sinh(\kappa(\omega)(X\sb{2}-X\sb{1}))/\kappa(\omega)
$
could only vanish at the points
$\omega=\pm\omega\sb{1,n}$,
where
\[
\omega\sb{J,n}
:=\sqrt{\frac{\pi^2 n^2}
{\abs{X\sb{J+1}-X\sb{J}}^2}+m^2},
\qquad
1\le J\le N-1,
\quad
n\in\N.
\]
Due to Assumption~\ref{ass-delta-small} and Lemma~\ref{lemma-m-m},
$\supp\hat\beta'(X\sb{1}+0,\omega)\cap\{\pm\omega\sb{1,n}\sothat n\in\N\}
=\emptyset$.
Hence, the middle line of (\ref{beta-beta-c})
at $x=X\sb{2}-0$
and the inequalities (\ref{imp})
imply that
\begin{equation}\label{sigma2-sigma1}
\inf\Sigma\sb{2}
= \inf\supp\hat g\sb{1}
<\inf\Sigma\sb{1},
\qquad
\sup\Sigma\sb{2}
= \sup\supp\hat g\sb{1}>\sup\Sigma\sb{1}.
\end{equation}
We proceed by induction, proving that
\begin{equation}\label{induction}
\inf\Sigma\sb{1}>\inf\Sigma\sb{2}>\dots>\inf\Sigma\sb{N},
\qquad
\sup\Sigma\sb{1}<\sup\Sigma\sb{2}<\dots<\sup\Sigma\sb{N}.
\end{equation}
It then follows that $\inf\Sigma\sb{N}<\sup\Sigma\sb{N}$.
Starting from $J=N$ and going to the left,
we also prove the opposite inequalities:
\begin{equation}\label{induction-2}
\inf\Sigma\sb{1}<\inf\Sigma\sb{2}<\dots<\inf\Sigma\sb{N},
\qquad
\sup\Sigma\sb{1}>\sup\Sigma\sb{2}>\dots>\sup\Sigma\sb{N}.
\end{equation}
The contradiction of (\ref{induction}) and (\ref{induction-2})
shows that
our assumption that $\inf\Sigma\sb{1}<\sup\Sigma\sb{1}$
was false,
hence
$\Sigma\sb{1}=\{\omega\sb{+}\}$
for some $\omega\sb{+}\in[-m,m]$.
\end{proof}

Thus, 
$\supp\hat\upbeta\sb{1}(\omega)=\Sigma\sb{1}\subset\{\omega\sb{+}\}$,
with $\omega\sb{+}\in[-m,m]$.
Therefore,
\begin{equation}\label{delom}
\hat\upbeta\sb{1}(\omega)=a\sb 1\delta(\omega-\omega\sb{+}),
\qquad {\rm with\ some}\ \ a\sb 1\in\C.
\end{equation}
Note that
the derivatives $\delta\sp{(k)}(\omega-\omega\sb{+})$, $k\ge 1$
do not enter the expression for
$\hat\upbeta\sb{1}(\omega)=\mathcal{F}\sb{t\to\omega}[\beta(X\sb 1,t)]$
since $\beta(x,t)$
is a bounded continuous function of $(x,t)\in\R^2$
due to the bound (\ref{beta-beta}).

\begin{lemma}\label{lem}
$\hat\beta(x,\omega)=a(x)\delta(\omega-\omega\sb{+})$, 
where $a(x)$ is a bounded continuous function.
\end{lemma}

\begin{proof}
For $x\le X\sb{1}$, the representation stated in the lemma 
follows from the first line in (\ref{beta-beta-c})
and from (\ref{delom}).
Let us prove this representation for $X\sb{1}\le x\le X\sb{2}$.
By (\ref{delom}), 
we have $\upbeta\sb{1}(t):=\beta(X\sb 1,t)=a\sb{1}e^{-i\omega\sb{+}t}/2\pi$,
hence
$g\sb{1}(t):=F\sb{1}(\upbeta\sb{1}(t))=b\sb{1}e^{-i\omega\sb{+}t}$
for some $b\sb{1}\in\C$
due to the $U(1)$-invariance (\ref{inv-f}).
Therefore, $\hat g\sb{1}(\omega)=2\pi b\sb{1}\delta(\omega-\omega\sb{+})$.
Moreover, by (\ref{beta-beta-c}), we have
$\hat\beta'(X\sb{1}-0,\omega)
=\kappa(\omega\sb{+})a\sb 1\delta(\omega-\omega\sb{+})$.
Hence,
the jump condition (\ref{jump-condition}) implies that
$\hat\beta'(X\sb{1}+0,\omega)=c\sb{1}\delta(\omega-\omega\sb{+})$,
for some $c\sb 1\in\C$.
Finally,  (\ref{beta-beta-c}) implies that
$\hat\beta(x,\omega)=a(x)\delta(\omega-\omega\sb{+})$
for $x\in [X\sb{1},X\sb{2}]$,
with $a(x)$ a continuous complex-valued function of $x$.
Proceeding by induction, we obtain 
similar representation for $\hat\beta(x,\omega)$ for all $x\in\R$.
\end{proof}

Now we can finish the proof of Proposition~\ref{prop-one-omega}.
Lemma \ref{lem} implies that
$\beta(x,t)=\phi(x)e^{-i\omega\sb{+}t}$,
where
$\phi(x)=a(x)/2\pi$.
We conclude
from (\ref{beta-beta})
that $\phi\in H\sp{1}(\R)$,
finishing the proof of Proposition~\ref{prop-one-omega}.
Note that
$\omega=\pm m$ could only correspond to the zero solution
(see Remark~\ref{remark-zero}).
\end{proof}

According to Remark~\ref{remark-end},
Proposition~\ref{prop-one-omega}
completes the proof of Theorem~\ref{main-theorem}.

\section{Multifrequency solitary waves}
\label{sect-counterexamples}
We will show that when the assumptions of Theorem~\ref{main-theorem}
are not satisfied, then the attractor could be more complicated
because the equation admits multifrequency solitary wave solutions.

\subsection{Wide gaps}
\label{sect-example-w}

Let us consider equation (\ref{kg-no}) with $N=2$,
under Assumptions~\ref{ass-a} and ~\ref{ass-nonl}.

\begin{proposition}
If the Assumption~\ref{ass-delta-small}
is violated, then the conclusion of Theorem~\ref{main-theorem}
may no longer be correct.
\end{proposition}

\begin{proof}
We will show that
if $L:=X\sb 2-X\sb 1$ is sufficiently large,
then
one can take $F\sb 1(\psi)$ and $F\sb 2(\psi)$
satisfying
Assumptions~\ref{ass-a} and~\ref{ass-nonl}
such that
the global attractor of the equation
contains the multifrequency solutions
which do not converge
to solitary waves of the form (\ref{solitary-waves}).
For our convenience, we assume that $X\sb{1}=0$, $X\sb{2}=L$.
We consider the model (\ref{kg-no}) with the nonlinearity
\begin{equation}
F\sb{1}(\psi)=F\sb{2}(\psi)
=F(\psi),
\qquad
{\rm where}
\quad
F(\psi)
=\alpha\psi+\beta\abs{\psi}^2\psi,
\qquad \alpha,\,\beta\in\R.
\end{equation}
In terms of the condition (\ref{f-is-such}), $p\sb{1}=p\sb{2}=2$.
We take $L$ to be large enough:
\begin{equation}\label{m-pi}
L>\frac{\pi}{2\sp{3/2}m}.
\end{equation}
Consider the function
\begin{equation}\label{psi-plus}
\psi(x,t)=
A (e^{-\kappa(\omega)\abs{x}}
+
e^{-\kappa(\omega)\abs{x-L}})
\sin(\omega t)
+B\chi\sb{[0,L]}(x)\sin(k(3\omega) x)\,\sin(3\omega t),
\qquad
A,\,B\in\C.
\end{equation}
Then
$\psi(x,t)$ solves (\ref{kg-no}) for $x$
away from the points $X\sb{J}$.
We require that
\begin{equation}\label{k-d}
k(3\omega)=\frac{\pi }{L},
\end{equation}
so that $\psi(x,t)$ is continuous in $x\in\R$
and symmetric with respect to $x=L/2$:
\[
\psi(x,t)
=
\psi(\frac{L}{2}-x,t),
\qquad x\in\R.
\]
We need
$\abs{\omega}<m$ to have $\kappa(\omega)>0$, and
$3\abs{\omega}>m$ to have $k(3\omega)\in\R$.
We take $\omega>0$, and thus
$m<3\omega<3m$.
By (\ref{k-d}), this means that we need
\[
m<\sqrt{\frac{\pi^2 }{L^2}+m^2}<3 m.
\]
The second inequality is satisfied
by (\ref{m-pi}).

Due to the symmetry of $\psi(x,t)$
with respect to $x=L/2$,
the jump condition
(\ref{jump-condition})
both at $x=0$ and at $x=L$ takes
the following identical form:
\begin{equation}\label{j-at-0-l}
2A\kappa(\omega) \sin\omega t-B k(3\omega)\sin 3\omega t
=
F\big(A (1+e^{-\kappa(\omega)L})\sin(\omega t)\big).
\end{equation}
Using the identity
\begin{equation}\label{trig}
\sin^3\theta=\frac 3 4\sin\theta-\frac 1 4\sin 3\theta,
\end{equation}
we see that
\begin{equation}
F(A (1+e^{-\kappa(\omega)L})\sin\omega t)
=
\Big(\alpha A(1+e^{-\kappa(\omega)L})
+
\frac{3}{4}
\beta \abs{A}^2 A(1+e^{-\kappa(\omega)L})^3
\Big)
\sin(\omega t)
-\frac{1}{4}\beta \abs{A}^2 A(1+e^{-\kappa(\omega)L})^3
\sin(3\omega t).
\end{equation}
Collecting in (\ref{j-at-0-l})
the terms at $\sin\omega t$ and at $\sin 3\omega t$,
we obtain the following system:
\begin{equation}\label{c-e}
\left\{
\begin{array}{l}
2A\kappa(\omega)
=\alpha A(1+e^{-\kappa(\omega)L})
+\frac 3 4\beta \abs{A}^2 A(1+e^{-\kappa(\omega)L})^3,
\\
B k(3\omega)=\frac 1 4\beta \abs{A}^2 A(1+e^{-\kappa(\omega)L})^3.
\end{array}
\right.
\end{equation}
Assuming that $A\ne 0$,
we divide the first equation by $A$:
\begin{equation}
2\kappa(\omega)
=\alpha (1+e^{-\kappa(\omega)L})
+\frac 3 4\beta \abs{A}^2(1+e^{-\kappa(\omega)L})^3.
\end{equation}
The condition for the existence of a solution $A\ne 0$
is
\begin{equation}\label{2k-alpha}
\Big(\frac{2\kappa(\omega)}{1+e^{-\kappa(\omega)L}}-\alpha\Big)\beta>0.
\end{equation}
Once we found $A$,
the second equation
in (\ref{c-e})
can be used to express $B$ in terms of $A$.

\begin{remark}
Condition (\ref{2k-alpha}) shows that we can choose
$\beta<0$ taking large $\alpha>0$.
The corresponding potential
$U(\psi)=-\alpha\abs{\psi}^2/2-\beta\abs{\psi}^4/4$
satisfies (\ref{bound-below})
and Assumptions~\ref{ass-a} and~\ref{ass-nonl}.
\end{remark}

\end{proof}

\subsection{Linear degeneration}
\label{sect-example-li}
Let us consider equation (\ref{kg-no})
with $N=2$,
under Assumptions~\ref{ass-a} and ~\ref{ass-delta-small}.

\begin{proposition}
If the Assumption~\ref{ass-nonl}
is violated, then the conclusion of Theorem~\ref{main-theorem}
may no longer be correct.
\end{proposition}

\begin{proof}
Again, we construct
multifrequency solutions.
Consider the equation
\begin{equation}
\ddot\psi
=\psi''-m^2\psi
+\delta(x)F\sb{1}(\psi)
+\delta(x-L)F\sb{2}(\psi),
\end{equation}
where
\begin{equation}
F\sb{1}(\psi)=\alpha\psi+\beta\abs{\psi}^2\psi,
\qquad
F\sb{2}(\psi)=\gamma\psi,
\qquad
\alpha,\ \beta,\ \gamma\in\R.
\end{equation}
Note that the function $F\sb{2}$ is linear,
failing to satisfy
Assumption~\ref{ass-nonl}.
The function
\[
\psi(x,t)
=
\left\{
\begin{array}{l}
(A+B) e^{\kappa(\omega)x}\sin(\omega t),
\qquad x\le 0,
\\
\big(A e^{-\kappa(\omega)x}+B e^{\kappa(\omega)x}\big)
\sin(\omega t)
+C\sinh(\kappa(3\omega)x)\sin(3\omega t),
\qquad x\in[0,L],
\\
(A e^{-\kappa(\omega)}+B e^{\kappa(\omega)(2L-x)})\sin(\omega t)
+\frac{C}{\sinh(\kappa(3\omega)L)}
e^{-\kappa(3\omega)(x-L)}\sin(3\omega t)
,
\qquad
x\ge L,
\end{array}
\right.
\]
where $\omega\in(0,m/3)$,
will be a solution if
the jump conditions are satisfied
at $x=0$ and at $x=L$:
\begin{equation}\label{j-at-0}
-\psi'(0+,t)+\psi'(0-,t)=\alpha\psi(0,t)+\beta\psi^3(0,t),
\end{equation}
\begin{equation}\label{j-at-l}
-\psi'(L+,t)+\psi'(L-,t)=\alpha\psi(L,t)+\beta\psi^3(L,t).
\end{equation}
We use the identity
\[
\alpha(A+B)\sin(\omega t)
+\beta((A+B)\sin(\omega t))^3
=\Big(\alpha(A+B)+\beta\frac{3(A+B)^3}{4}\Big)\sin(\omega t)
-\beta\frac{(A+B)^3}{4}\sin(3\omega t)
\]
which follows from (\ref{trig}).
Collecting the terms at $\sin(\omega t)$
and at $\sin(3\omega t)$,
we write the condition (\ref{j-at-0})
as the following system of equations:
\begin{eqnarray}
&&
2\kappa(\omega)A=\Big(\alpha(A+B)+\beta\frac{3(A+B)^3}{4}\Big),
\label{c01}
\\
&&
-\kappa(3\omega)C=-\beta\frac{(A+B)^3}{4}.
\label{c03}
\end{eqnarray}
Similarly, the condition (\ref{j-at-l}) is equivalent to the
following two equations:
\begin{eqnarray}
&&
2B\kappa(\omega)e^{\kappa(\omega)L}
=\gamma(A e^{-\kappa(\omega)L}+B e^{\kappa(\omega)L}),
\label{cl1}
\\
&&
\frac{\kappa(3\omega)C}{\sinh(\kappa(3\omega)L)}
+
\kappa(3\omega)C\cosh(\kappa(3\omega)L)
=\gamma C\sinh(\kappa(3\omega)L).
\label{cl3}
\end{eqnarray}
Equations (\ref{c01}), (\ref{c03}), (\ref{cl1}), and (\ref{cl3})
could be satisfied for arbitrary $L>0$.
Namely, for any $\omega\in(0,m/3)$,
one uses (\ref{cl3}) to determine $\gamma$.
For any $\beta\ne 0$,
there is always a solution $A$, and $B$
to the nonlinear system (\ref{c01}), (\ref{cl1}).
Finally, $C$ is obtained from (\ref{c03}).

\end{proof}

\appendix

\section{Global well-posedness}

\label{sect-existence}

Here we prove Theorem~\ref{theorem-well-posedness}.
We first need
to adjust the nonlinearity $F$
so that it becomes bounded, together with
its derivatives.
Define
\begin{equation}\label{def-lambda-0}
\lambda\sb{0}
=\sqrt{\frac{\mathcal{H}(\psi\sb{0},\pi\sb{0})-\sum\sb{J}{A}\sb{J}}{m-\sum\sb{J}{B}\sb{J}}},
\end{equation}
where $(\psi\sb{0},\pi\sb{0})\in{\E}$ is the initial data
from Theorem~\ref{theorem-well-posedness}
and ${A}\sb{J}$, ${B}\sb{J}$ are constants from (\ref{bound-below}).
Then we may pick a modified potential function
$\widetilde{U}\sb{J}\in C\sp 2(\C,\R)$,
$\widetilde{U}\sb{J}(\psi)=\widetilde{U}\sb{J}(\abs{\psi})$,
$j=1,\,2$,
so that
\begin{equation}\label{new-U}
\widetilde{U}\sb{J}(\psi)=U\sb{J}(\psi)
\qquad{\rm for}\ \abs{\psi}\le\lambda\sb{0},
\quad
\psi\in\C,
\end{equation}
$\widetilde{U}\sb{J}(\psi)$ satisfy (\ref{bound-below})
with the same constants ${A}\sb{J}$, ${B}\sb{J}$ as $U\sb{J}(\psi)$ do:
\begin{equation}\label{new-U-2}
\widetilde{U}\sb{J}(\psi)\ge {A}\sb{J}-{B}\sb{J}\abs{\psi}^2,
\quad{\rm for}\ \psi\in\C,\quad
{\rm where}
\quad
{A}\sb{J}\in\R,
\quad
{B}\sb{J}\ge 0,
\quad
1\le J\le N,
\quad
\sum\sb{J}{B}\sb{J}<m,
\end{equation}
and so that
$\abs{\widetilde{U}\sb{J}(\psi)}$,
$\abs{\widetilde{U}\sb{J}'(\psi)}$,
and $\abs{\widetilde{U}\sb{J}''(\psi)}$
are bounded for $\psi\ge 0$.
We define
\begin{equation}\label{f-reg}
\widetilde{F}\sb{J}(\psi)
=-\nabla \widetilde{U}\sb{J}(\psi),
\qquad\psi\in\C,
\end{equation}
where $\nabla$ denotes the gradient with respect to
$\Re\psi$, $\Im\psi$;
Then
$\widetilde{F}\sb{J}(e\sp{is}\psi)
=e\sp{is}\widetilde{F}\sb{J}(\psi)$ for any $\psi\in\C$, $s\in\R$.

We consider the Cauchy problem of type
(\ref{kg-no}) with the modified
nonlinearity,
\begin{equation}\label{kg-no-a}
\left\{
\begin{array}{l}
\ddot\psi(x,t)
=\psi''(x,t)-m^2\psi(x,t)+\sum\sb{J}\delta(x-X\sb{J})
\widetilde F\sb{J}(\psi(X\sb{J},t)),
\qquad
x\in\R,
\quad t\in\R,
\\
\psi\at{t=0}=\psi\sb{0}(x),
\qquad
\dot\psi\at{t=0}=\pi\sb{0}(x).
\end{array}\right.
\end{equation}
Equation (\ref{kg-no-a})
formally can be written as the following Hamiltonian system
(Cf. (\ref{kg-no-h})):
\begin{equation}\label{kg-no-h-t}
\dot\Psi(t)=\mathcal{J}\,D\widetilde{\mathcal{H}}(\Psi),
\qquad
\mathcal{J}=\left[\begin{array}{cc}0&1\\-1&0\end{array}\right],
\end{equation}
where $D\widetilde{\mathcal{H}}$ is the variational
derivative of the Hamilton functional
\begin{equation}\label{kg-no-a-h}
\widetilde{\mathcal{H}}(\Psi)
=\int\limits\sb{\R}
\left(
\abs{\pi}\sp 2
+\abs{\nabla\psi}\sp 2+m^2\abs{\psi}\sp 2
\right)\,dx
+\sum\sb{J}\widetilde{U}\sb{J}(\psi(X\sb{J},t)),
\quad
\Psi=\left[\begin{array}{c}\psi(x)\\
\pi(x)\end{array}\right]\in\E,
\end{equation}
which is Fr\'echet differentiable in the space
${\E}=H\sp{1}\times L\sp 2$.
By the Sobolev embedding theorem,
$L\sp\infty(\R)\subset H\sp{1}(\R)$,
and there is the following inequality:
\begin{equation}\label{sobolev-embedding}
\norm{\psi}\sb{L\sp\infty}^2
\le\frac{1}{2m}
(\norm{\psi'}\sb{L\sp 2}^2+m^2\norm{\psi}\sb{L\sp 2}^2)
\le\frac{1}{2m}
\norm{\Psi}\sb{\E}^2.
\end{equation}
Thus, (\ref{new-U-2}) leads to
\begin{equation}\label{bound-on-u}
\widetilde{U}\sb{J}(\psi(0))
\ge {A}\sb{J}-{B}\sb{J}\norm{\psi}\sb{L\sp\infty}^2
\ge {A}\sb{J}-\frac{{B}\sb{J}}{2m}\norm{\Psi}\sb{\E}^2.
\end{equation}
Taking into account (\ref{kg-no-a-h}),
we obtain the inequality
\begin{equation}
\norm{\Psi}\sb{\E}\sp 2
=2\widetilde{\mathcal{H}}(\Psi)-2\sum\sb{J}\widetilde{U}\sb{J}(\psi(X\sb{J}))
\le 2\widetilde{\mathcal{H}}(\Psi)-2\sum\sb{J}{A}\sb{J}+\frac{\sum\sb{J}{B}\sb{J}}{m}
\norm{\Psi}\sb{\E}^2,
\qquad
\Psi\in\E.
\end{equation}
It follows that
\begin{equation}\label{t-bound-1}
\norm{\Psi}\sb{\E}\sp 2
\le\frac{2m}{m-\sum\sb{J}{B}\sb{J}}
\Big(\widetilde{\mathcal{H}}(\Psi)-\sum\sb{J}{A}\sb{J}\Big),
\qquad\Psi\in\E.
\end{equation}

\begin{lemma}\label{lemma-same-u}
\begin{enumerate}
\item
There is the identity
$\widetilde{\mathcal{H}}(\Psi\sb{0})
=\mathcal{H}(\Psi\sb{0})$.
\item
If $\Psi=\left[\begin{array}{c}\psi(x)\\\pi(x)\end{array}\right]\in{\E}$
satisfies
$\widetilde{\mathcal{H}}(\Psi)\le\widetilde{\mathcal{H}}(\Psi\sb{0})$,
then
$\ \widetilde{U}\sb{J}(\psi(x))=U\sb{J}(\psi(x))$
for any $x\in\R$.
\end{enumerate}
\end{lemma}

\begin{proof}
According to (\ref{t-bound-1}),
the Sobolev embedding (\ref{sobolev-embedding}),
and the choice of $\lambda\sb{0}$ in (\ref{def-lambda-0}),
\begin{equation}
\norm{\psi\sb{0}}\sb{L\sp\infty}^2
\le\frac{1}{2m}\norm{\Psi\sb{0}}\sb{\E}^2
\le\frac{\mathcal{H}(\Psi\sb{0})-\sum\sb{J}{A}\sb{J}}{m-\sum\sb{J}{B}\sb{J}}
=\lambda\sb{0}^2.
\end{equation}
Thus, by (\ref{new-U}),
$
\widetilde{U}(\psi\sb{0}(x))=U(\psi\sb{0}(x))
$
for all $x\in\R$.
This proves ({\it i}).

By (\ref{sobolev-embedding}),
the relation (\ref{t-bound-1}),
the condition
$\widetilde{\mathcal{H}}(\Psi)\le\widetilde{\mathcal{H}}(\Psi\sb{0})$,
and part ({\it i}) of the Lemma, we have:
\[
\norm{\psi}\sb{L\sp\infty}^2
\le\frac{1}{2m}\norm{\Psi}\sb{\E}^2
\le\frac{\widetilde{\mathcal{H}}(\Psi)-\sum\sb{J}{A}\sb{J}}{m-\sum\sb{J}{B}\sb{J}}
\le\frac{\widetilde{\mathcal{H}}(\Psi\sb{0})-\sum\sb{J}{A}\sb{J}}{m-\sum{B}\sb{J}}
=\frac{\mathcal{H}(\Psi\sb{0})-\sum\sb{J}{A}\sb{J}}{m-\sum\sb{J}{B}\sb{J}}
=\lambda\sb{0}^2.
\]
Now the statement ({\it ii}) follows by (\ref{new-U}).
\end{proof}

If $\Psi(t)$ solves (\ref{kg-no-h-t}),
then
$\widetilde{\mathcal{H}}(\Psi(t))=\widetilde{\mathcal{H}}(\Psi\sb{0})$,
By Lemma~\ref{lemma-same-u}~({\it ii}),
$\widetilde{U}\sb{J}(\psi(x,t))=U\sb{J}(\psi(x,t))$
for all $x\in\R$, $t\in\R$.
Hence, $\widetilde{F}\sb{J}(\psi(x,t))=F\sb{J}(\psi(x,t))$
for all $x\in\R$, $t\ge 0$,
allowing us to conclude that
$\psi(t)$ solves (\ref{kg-no})
as well as (\ref{kg-no-a}).
The rest of the proof
of Theorem~\ref{theorem-well-posedness}
repeats the proof of a similar result for the case $N=1$
\cite[Theorem 2.3]{ubk-arma}.

\bibliographystyle{sima}
\bibliography{comech,ubk-mathsci,ubk-local,all}

\end{document}